\documentclass[12pt]{article}
\usepackage{amsmath,amsfonts,epsfig}
\title{\bf Chordal Coxeter Groups}
\author{John Ratcliffe and Steven Tschantz \\ 
Mathematics Department, Vanderbilt University, \\
Nashville TN 37240, USA}
\newtheorem{theorem}{Theorem}[section]
\newtheorem{proposition}[theorem]{Proposition}
\newtheorem{lemma}[theorem]{Lemma}

\newtheorem{conjecture}[theorem]{Conjecture}

\newenvironment{proof}{{\bf Proof:\ }}{\hfill$\square$\vspace{.2in}}

\def\ov{\overline}

\date{}
\begin{document}
\maketitle

\noindent {\bf Abstract:} A solution of the isomorphism problem is presented 
for the class of Coxeter groups $W$ that have  
a finite set of Coxeter generators $S$ such that the underlying graph 
of the presentation diagram of the system $(W,S)$ has the property that
every cycle of length at least four has a cord. 
As an application, we construct counterexamples to two main conjectures 
concerning the isomorphism problem for Coxeter groups. 

\section{Introduction} 

The isomorphism problem for finitely generated Coxeter groups is the problem of deciding 
if two finite Coxeter matrices define 
isomorphic Coxeter groups. Coxeter \cite{Coxeter} solved this problem for finite irreducible 
Coxeter groups.  
Recently there has been considerable interest and activity on the isomorphism problem 
for arbitrary finitely generated Coxeter groups. 
For a recent survey, see M\"uhlherr \cite{Muhlherr}. 

A graph is said to be {\it chordal} if every cycle of length at least four has a cord. 
The class of finite chordal graphs is the smallest class of finite graphs 
that contains all complete finite graphs and is closed with respect 
to taking the union of two graphs intersecting along a complete graph. 
The class of finite chordal graphs is an important class of graphs 
that has been extensively studied and is well understood. 
As a reference for chordal graphs, 
see Brandst\"adt, Le, and Spinrad \cite{B-L-S}. 

A Coxeter group $W$ is said to be {\it chordal} if $W$ has 
a set of Coxeter generators $S$ such that the underlying graph 
of the presentation diagram of the system $(W,S)$ is chordal. 
In this paper, we present a solution of the isomorphism problem 
for the class of finitely generated chordal Coxeter groups. 
As an application, we construct a counterexample to Conjecture 8.1 of 
Brady, McCammond, M\"uhlherr, and Neumann \cite{B-M-M-N} and to Conjecture 1 in 
M\"uhlherr's survey \cite{Muhlherr} by employing a new type of diagram twisting.  
On the positive side, we prove Conjecture 2 in M\"uhlherr's survey 
for the class of finitely generated chordal Coxeter groups. 

In \S 2, we state some preliminary results. 
In \S 3, we discuss twisting Coxeter systems. 
In \S 4, we prove the automorphic angle deformation theorem. 
In \S 5, we prove the star decomposition theorem. 
In \S 6, we prove the angle deformation theorem. 
In \S 7, we prove some lemmas about conjugating visual subgroups. 
In \S 8, we prove the sharp angle theorem. 
In \S 9, we prove the twist equivalence theorem.  
In \S 10, we present our solution of the isomorphism problem for finitely generated 
chordal Coxeter groups. 
In \S 11, we present our counterexample to Conjecture 8.1 in \cite{B-M-M-N} and 
Conjecture 1 in M\"uhlherr's survey.

\section{Preliminaries} 
 
A {\it Coxeter matrix} is a symmetric matrix $M = (m(s,t))_{s,t\in S}$ 
with $m(s,t)$ either a positive integer or infinity and $m(s,t) =1$ 
if and only if $s=t$. A {\it Coxeter system} with Coxeter matrix $M = (m(s,t))_{s,t\in S}$ 
is a pair $(W,S)$ consisting of a group $W$ and a set of generators $S$ for $W$ 
such that $W$ has the presentation
$$W =\langle S \ |\ (st)^{m(s,t)}:\, s,t \in S\ \hbox{and}\ m(s,t)<\infty\rangle$$
If $(W,S)$ is a Coxeter system with Coxeter matrix $M = (m(s,t))_{s,t\in S}$, 
then the order of $st$ is $m(s,t)$ for each $s,t$ in $S$ by 
Prop. 4, p. 92 of Bourbaki \cite{Bourbaki}, and so a Coxeter system $(W,S)$ 
determines its Coxeter matrix; moreover, any Coxeter matrix 
$M = (m(s,t))_{s,t\in S}$ determines a Coxeter system $(W,S)$ where $W$ 
is defined by the above presentation. If $(W,S)$ is a Coxeter system, 
then $W$ is called a {\it Coxeter group} and $S$ is called a set of {\it Coxeter generators} 
for $W$, and the cardinality of $S$ is called the {\it rank} of $(W,S)$. 
It follows from Theorem 2 (i), p. 20 of Bourbaki  \cite{Bourbaki} 
that a Coxeter system $(W,S)$ has finite rank if and only if $W$ 
is finitely generated.

Let $(W,S)$ be a Coxeter system. A {\it visual subgroup} of $(W,S)$ 
is a subgroup of $W$ of the form $\langle A\rangle$ for some $A \subset S$. 
If $\langle A\rangle$ is a visual subgroup of $(W,S)$, 
then $(\langle A\rangle, A)$ is also a Coxeter system 
by Theorem 2 (i), p. 20 of Bourbaki  \cite{Bourbaki}.

When studying a Coxeter system $(W,S)$ with Coxeter matrix $M$ 
it is helpful to have a visual representation of $(W,S)$. 
There are two graphical ways of representing $(W,S)$ 
and we shall use both depending on our needs. 

The {\it Coxeter diagram} ({\it {\rm C}-diagram}) of $(W,S)$ is the labeled undirected graph 
$\Delta = \Delta(W,S)$ with vertices $S$ and edges 
$$\{(s,t) : s, t \in S\ \hbox{and}\ m(s,t) > 2\}$$
such that an edge $(s,t)$ is labeled by $m(s,t)$. 
Coxeter diagrams are useful for visually representing finite Coxeter groups. 
If $A\subset S$, then $\Delta(\langle A\rangle,A)$ is the 
subdiagram of $\Delta(W,S)$ induced by $A$.

A Coxeter system $(W,S)$ is said to be {\it irreducible} 
if its C-diagram $\Delta$ is connected. 
A visual subgroup $\langle A\rangle$ of $(W,S)$ is said to be {\it irreducible} 
if $(\langle A\rangle, A)$ is irreducible. 
A subset $A$ of $S$ is said to be {\it irreducible} if $\langle A\rangle$ is irreducible. 

A subset $A$ of $S$ is said to be a {\it component} of $S$ if $A$ is a maximal irreducible 
subset of $S$ or equivalently if $\Delta(\langle A\rangle, A)$ is a connected component 
of $\Delta(W,S)$. 
The connected components of the $\Delta(W,S)$
represent the factors of a direct product decomposition of $W$.

The {\it presentation diagram} ({\it {\rm P}-diagram}) of $(W,S)$ is the labeled undirected graph 
$\Gamma = \Gamma(W,S)$ with vertices $S$ and edges 
$$\{(s,t) : s, t \in S\ \hbox{and}\ m(s,t) < \infty\}$$
such that an edge $(s,t)$ is labeled by $m(s,t)$. 
Presentation diagrams are useful for visually representing infinite Coxeter groups. 
If $A\subset S$, then $\Gamma(\langle A\rangle,A)$ is the 
subdiagram of $\Gamma(W,S)$ induced by $A$. 
The connected components of $\Gamma(W,S)$
represent the factors of a free product decomposition of $W$.

Let $(W,S)$ and $(W',S')$ be Coxeter systems 
with P-diagrams $\Gamma$ and $\Gamma'$, respectively. 
An {\it isomorphism} $\phi: (W,S) \to (W',S')$ of Coxeter systems 
is an isomorphism $\phi: W\to W'$ such that $\phi(S) = S'$. 
An {\it isomorphism} $\psi:\Gamma\to \Gamma'$ of P-diagrams is a bijection 
from $S$ to $S'$ that preserves edges and their labels. 
\begin{proposition} 
Let $(W,S)$ and $(W',S')$ be Coxeter systems with P-dia\-grams $\Gamma$ and $\Gamma'$, 
respectively. Then
\begin{enumerate}
\item $(W,S) \cong (W',S')$ if and only if $\Gamma \cong \Gamma'$,
\item $W\cong W'$ if and only if $W$ has a set of Coxeter generators $S''$ such that 
		$(W,S'')\cong (W',S')$,
\item $W\cong W'$ if and only if $W$ has a P-diagram $\Gamma''$ such that 
		$\Gamma''\cong \Gamma'$.
\end{enumerate}
\end{proposition} 
\begin{proof} (1) If $\phi: (W,S) \to (W',S')$ is an isomorphism, then 
$\phi$ restricts to an isomorphism $\overline\phi:\Gamma \to \Gamma'$ and 
if $\psi:\Gamma \to \Gamma'$ is an isomorphism, then 
$\psi$ extends to a unique isomorphism $\hat \psi: (W,S) \to (W',S')$. 

(2) If $\phi:W\to W'$ is an isomorphism, then $S''=\phi^{-1}(S')$ 
is a set of Coxeter generators for $W$ and $\phi:(W,S'')\to (W',S')$ 
is an isomorphism. 

(3) Statement (3) follows from (1) and (2).
\end{proof}

A Coxeter group $W$ is said to be {\it rigid} 
if for any two sets of Coxeter generators $S$ and $S'$ for $W$, 
there is an automorphism $\alpha: (W,S)\to (W,S')$ 
or equivalently any two sets of Coxeter generators $S$ and $S'$ for $W$ 
determine isomorphic P-diagrams for $W$. 
A Coxeter group $W$ is said to be {\it strongly rigid} 
if any two sets of Coxeter generators for $W$ are conjugate. 

A Coxeter system $(W,S)$ is said to be {\it complete} if 
the underlying graph of the P-diagram of $(W,S)$ is complete. 
A Coxeter system $(W,S)$ is said to be {\it finite} (resp.  
{\it infinite}) if $W$ is finite (resp. infinite).

\begin{theorem}{\rm (Caprace, Franzsen, 
Howlett, and M\"uhlherr \cite{C-M},\cite{F-H-M})} 
If $(W,S)$ is an infinite, complete, irreducible Coxeter system of finite rank, 
then $W$ is strongly rigid. 
\end{theorem}

We shall use Coxeter's notation on p. 297 of \cite{Coxeterb}
for the irreducible spherical Coxeter simplex reflection groups except 
that we denote the dihedral group ${\bf D}_2^k$ by ${\bf D}_2(k)$. 
Subscripts denote the rank of a Coxeter system in Coxeter's notation. 
Coxeter's notation partly agrees with but differs from Bourbaki's notation on p.193 of 
\cite{Bourbaki}.

Coxeter \cite{Coxeter} proved that every finite irreducible Coxeter system 
is isomorphic to exactly one 
of the Coxeter systems ${\bf A}_n$, $n\geq 1$, ${\bf B}_n$, $n\geq 4$, ${\bf C}_n$, 
$n\geq 2$, ${\bf D}_2(k)$, $k\geq 5$, 
${\bf E}_6$, ${\bf E}_7$, ${\bf E}_8$, ${\bf F}_4$, ${\bf G}_3$, ${\bf G}_4$. 
For notational convenience, we define ${\bf B}_3 = {\bf A}_3$, ${\bf D}_2(3) = {\bf A}_2$, 
and ${\bf D}_2(4) = {\bf C}_2$

The {\it type} of a finite irreducible Coxeter system $(W,S)$ is the 
isomorphism type of $(W,S)$ represented by one of the systems 
${\bf A}_n$, ${\bf B}_n$, ${\bf C}_n$, ${\bf D}_2(k)$, 
${\bf E}_6$, ${\bf E}_7$, ${\bf E}_8$, ${\bf F}_4$, ${\bf G}_3$, ${\bf G}_4$. 
The {\it type} of an irreducible subset $A$ of $S$  
is the type of $(\langle A\rangle,A)$.  

The C-diagram of ${\bf A}_n$ is a linear diagram 
with $n$ vertices and all edge labels 3. 
The C-diagram of ${\bf B}_n$ is a Y-shaped diagram 
with $n$ vertices and all edge labels 3 and two short 
arms of consisting of single edges. 
The C-diagram of ${\bf C}_n$ is a linear diagram with $n$ vertices 
and all edge labels 3 except for the last edge labelled 4. 
The C-diagram of ${\bf D}_2(k)$ is a single edge with label $k$. 
The C-diagrams of ${\bf E}_6$, ${\bf E}_7, {\bf E}_8$ 
are star shaped with three arms and all edge labels 3. 
One arm has length one and another has length two. 
The C-diagram of ${\bf F}_4$ is a linear diagram 
with edge labels $3,4,3$ in that order.  
The C-diagram of ${\bf G}_3$ is a linear diagram 
with edge labels $3, 5$. 
The C-diagram of ${\bf G}_4$ is a linear diagram 
with edge labels $3,3,5$ in that order.

\section{Twisting Coxeter Systems}  

Let $(W,S)$ be a Coxeter system of finite rank.  
Suppose that $S_1,S_2\subset S$, with $S=S_1\cup S_2$, and $S_0=S_1\cap S_2$ 
are such that there is no defining relator of $W$ (no edge of
the P-diagram) between an element of $S_1-S_0$ and $S_2-S_0$. 
Then we can write $W$ as a visual amalgamated product
$W=\langle S_1\rangle*_{\langle S_0\rangle}\langle S_2\rangle$. 
We say that $S_0$ {\it separates} $S$ if
$S_1-S_0\neq\emptyset$ and $S_2-S_0\neq\emptyset$.  
The amalgamated product decomposition of $W$ will be
nontrivial if and only if $S_0$ separates $S$. 
If $S_0$ separates $S$, we call the triple $(S_1,S_0,S_2)$ 
a {\it separation} of $S$. 
Note that $S_0$ separates $S$ if and only if $S_0$ 
separates $\Gamma(W,S)$, that is, there are $s_1,s_2$ in $S-S_0$ 
such that every path in $\Gamma(W,S)$ from $s_1$ to $s_2$ must pass 
through $S_0$. 

Let $\ell \in \langle S_0\rangle$ such that 
$\ell S_0\ell^{-1}=S_0$. 
By Lemma 4.5 of \cite{M-R-T}, we have $S_0=S_{\bullet}\cup (S_0-S_{\bullet})$
where $S_{\bullet}$ generates a finite group, each
element of $S_{\bullet}$ commutes with each element of $S_0-S_{\bullet}$, 
and $\ell$ is the longest element of $\langle S_{\bullet}\rangle$.  
The triple $(S_1,\ell, S_2)$ determines an elementary
twist of $(W,S)$ giving a new Coxeter
generating set $S_*=S_1\cup \ell S_2\ell^{-1}$ of $W$. 
We call $S_1$ the {\it left} (or {\it invariant}) side of the separation of $S$
of the twist determined by $(S_1,\ell, S_2)$. 

In application, it is simpler to consider a more general kind of twisting.  
Suppose $S_0$ and $\bar S_0\subset S_2$ generate conjugate subgroups of
$\langle S_2\rangle$.  Suppose $d\in \langle S_2\rangle$ is such
that $d\bar S_0d^{-1}=S_0$.  Then $S_1\cap dS_2d^{-1}=S_0$, since
$$S_0\subset S_1\cap dS_2d^{-1}\subset S_1\cap\langle S_2\rangle = S_0.$$
A {\it generalized twist} (or simply {\it twist}) of $(W,S)$ 
in this situation gives a new Coxeter generating
set $S_*=S_1\cup dS_2d^{-1}$ of $W$ and a new visual 
amalgamated product decomposition 
$W=\langle S_1\rangle*_{\langle S_0\rangle}\langle d S_2d^{-1}\rangle$

Elementary and generalized twists can be easily understood in terms 
of their effects on P-diagrams. 
The P-diagram of $(W,S)$ is the union of the 
P-diagrams for $\langle S_1\rangle$ and $\langle S_2\rangle$
overlapping in the P-diagram for $\langle S_0\rangle$. 
The P-diagram for $(W,S_*)$ is obtained from the P-diagram of $(W,S)$ 
by twisting the P-diagram of $\langle S_2\rangle$, that is, 
removing the P-diagram for $\langle S_2\rangle$, replacing it by the isomorphic 
P-diagram of $\langle dS_2d^{-1}\rangle$, and attaching it to the P-diagram
for $\langle S_1\rangle$ along $S_0 = d\bar S_0d^{-1}$.   
If $S_0 = \emptyset$, we call the twist {\it degenerate}. 
A degenerate twist does not change the isomorphism type of the P-diagram. 
This includes the case where
$S_1=S_0=\emptyset$, $S_2=S$, giving $S_*=dSd^{-1}$ the
conjugation of $S$ by an arbitrary $d\in W$.  
Any nondegenerate generalized twist of a Coxeter system $(W,S)$ can be realized 
by a finite sequence of elementary twists. 
Two Coxeter systems $(W,S)$ and $(W,S')$ are said to be {\it twist equivalent} 
if there is a sequence of twists that transforms the system $(W,S)$ into 
the system $(W,S')$. 

If $\Lambda$ is a visual graph of groups decomposition of $(W,S)$,  
then the graph of $\Lambda$ is a tree, since the abelianization of $W$ is finite.  
A graph of groups decomposition is said to be {\it reduced} 
if no edge group is equal to an incident vertex group. 
The following theorem was proved in \cite{M-R-T}.

\begin{theorem} 
{\rm (The Decomposition Matching Theorem)} 
Suppose $(W,S)$ and $(W,S')$ are Coxeter systems for the same
finitely generated Coxeter group and $W$ has a nontrivial splitting as $A*_CB$. 
Then $(W,S)$ and $(W,S')$ are twist equivalent to Coxeter systems $(W,R)$
and $(W,R')$, respectively, such that there exists a nontrivial
reduced visual graph of groups decomposition $\Psi$ of $(W,R)$
and a nontrivial reduced visual graph of groups decomposition
$\Psi'$ of $(W,R')$ having the same graphs and the same vertex
and edge groups and all edge groups equal and a subgroup of a
conjugate of $C$.
\end{theorem}

A {\it simplex} $A$ of $(W,S)$ is a subset $A$ of $S$ such that 
$(\langle A\rangle, A)$ is a complete Coxeter system. 
A simplex $A$ of $(W,S)$ is said to be {\it spherical} 
if $\langle A\rangle$ is finite.  

A Coxeter systems $(W,S)$ is said to be {\it chordal} 
if the underlying graph of the P-diagram of $(W,S)$ is chordal. 

\begin{lemma} 
Let $(W,S)$ be a chordal Coxeter system. 
Let $(S_1,S_0,S_2)$ be a separation of $S$ such that $S_0$ is a simplex. 
Let $d$ be an element of $\langle S_2\rangle$ such that $d^{-1}S_0 d \subset S_2$, 
and let $S_\ast = S_1 \cup dS_2d^{-1}$. 
Then the twisted system $(W,S_\ast)$ is chordal. 
\end{lemma}
\begin{proof}
Let $W_i = \langle S_i\rangle$ for $i=0,1,2$. 
We have 
$$\Gamma(W,S_\ast) = \Gamma(W_1,S_1)\cup \Gamma(dW_2d^{-1},dS_2d^{-1})$$
and
$$\Gamma(W_0,S_0) = \Gamma(W_1,S_1)\cap \Gamma(dW_2d^{-1},dS_2d^{-1}).$$
The underlying graphs of $\Gamma(W_1,S_1)$ and $\Gamma(dW_2d^{-1},dS_2d^{-1})$ are chordal 
and the underlying graph of $\Gamma(W_0,S_0)$ is complete. 
Therefore the underlying graph of $\Gamma(W,S_\ast)$ is chordal by Theorem 2 of Dirac \cite{Dirac}. 
\end{proof}

\begin{theorem} 
Let $(W,S)$ be a chordal Coxeter system of finite rank, 
and let $S'$ be another set of Coxeter generators for $W$. 
Then $(W,S')$ is chordal. 
\end{theorem}
\begin{proof}
The proof is by induction on $|S|$. 
The theorem is clear if $|S| =1$, 
so assume $|S| > 1$ and the theorem 
is true for all chordal Coxeter systems 
of rank less than $|S|$. 
Assume first that $(W,S)$ is complete. 
Then $(W,S')$ is complete by Proposition 5.10 of \cite{M-R-T}, 
and so $(W,S')$ is chordal. 

Now assume that $(W,S)$ is incomplete. 
Then $\Gamma(W,S)$ has a cut set. 
Let $C \subset S$ be a minimal cut set of $\Gamma(W,S)$. 
By Theorem 1 of Dirac \cite{Dirac}, 
we have that $C$ is a simplex.  
Let $(A,C,B)$ be the separation of $S$ determined by $C$. 
Then $W = \langle A\rangle \ast_{\langle C\rangle} \langle B\rangle$ 
is a nontrivial splitting. 
By the Decomposition Matching Theorem, 
$(W,S)$ and $(W,S')$ are twist equivalent 
to Coxeter systems $(W,R)$ and $(W,R')$, respectively, 
such that there exists a nontrivial
reduced visual graph of groups decomposition $\Psi$ of $(W,R)$
and a nontrivial reduced visual graph of groups decomposition
$\Psi'$ of $(W,R')$ having the same graphs and the same vertex
and edge groups and all edge groups equal and a subgroup of a
conjugate of $\langle C \rangle$.
The Coxeter systems $(W,R)$ and $(W,S)$ are twist equivalent, 
and so $|R| = |S|$ and $(W,R)$ is chordal by Lemma 3.2. 

Let $\{(W_i,R_i)\}_{i=1}^k$ be the Coxeter systems of the vertex groups of $\Psi$, 
and let $(W_0,R_0)$ be the Coxeter system of the edge group of $\Psi$. 
Then $k\geq 2$, and $R =\cup_{i=1}^k R_i$, and $\cap_{i=1}^k R_i = R_0$, 
and $R_i - R_0 \neq \emptyset$ for each $i > 0$, 
and $m(a,b) = \infty$ for each $a \in R_i-R_0$ and $b \in R_j-R_0$ with $i\neq j$. 
From the proof of Lemma 8.4 of \cite{M-R-T}, we deduce that $R_0$ 
is conjugate to a subset of $S$. 
By Lemma 4.3 of \cite{M-R-T}, we have that $R_0$ is conjugate to a subset of $C$, 
and so $R_0$ is a simplex. 
Let $\{(W_i',R_i')\}_{i=1}^k$ be the Coxeter systems of the vertex groups of $\Psi'$ 
indexed so that $W_i' = W_i$ for each $i$, 
and let $(W_0',R_0')$ be the Coxeter system of the edge group of $\Psi'$. 
Then $W_0' = W_0$, and $R' =\cup_{i=1}^k R_i'$, and $\cap_{i=1}^k R_i' = R_0'$, 
and $R_i' - R_0' \neq \emptyset$ for each $i > 0$, 
and $m(a',b') = \infty$ for each $a' \in R_i'-R_0'$ and $b' \in R_j'-R_0'$ with $i\neq j$. 

The Coxeter system $(W_i,R_i)$ is chordal and $|R_i| < |R|$ for each $i$. 
Hence by the induction hypothesis, $(W_i,R_i')$ is chordal for each $i$. 
As $(W,R_0)$ is complete, $(W,R_0')$ is complete by Prop. 5.10 of \cite{M-R-T}. 
This implies that $(W,R')$ is chordal by Theorem 2 of Dirac \cite{Dirac} 
and induction on $k$. 
Therefore $(W,S')$ is chordal, since $(W,S')$ is twist equivalent to $(W,R')$. 
\end{proof}

\section{The Automorphic Deformation Theorem}

Let $(W,S)$ be a Coxeter system, 
and let $A \subset S$. Define 
$$A^\perp = \{s \in S \ |\ m(s,a) = 2\ \hbox{for all}\ a \in A\}.$$ 

Let $\{c,d\} \subset S$, and let $B \subset S$. 
We say that $B$ is a $(c,d)$-{\it separator} of $S$ 
if there is a separation $(S_1,B,S_2)$ of $S$ such that $c \in S_1-B$ 
and $d \in S_2 - B$. 
Let $\Gamma$ be the P-diagram of $(W,S)$. 
Then $B$ is a $(c,d)$-separator if and only if $c$ and $d$ lie in 
different connected components of $\Gamma-B$. 

\begin{lemma} 
Let $(W,S)$ be a chordal Coxeter system, and let $A = \{a,b,c\}$ 
be a maximal irreducible simplex of $(W,S)$, 
and let $d\in S-(A\cup A^\perp)$ such that $\{a,b,d\}$ is a simplex of $(W,S)$. 
Then $\{a,b\}\cup A^\perp$ is a $(c,d)$-separator. 
\end{lemma}
\begin{proof}
On the contrary, suppose $\{a,b\}\cup A^\perp$ is not a $(c,d)$-separator. 
Then there is a path in $\Gamma(W,S)$ from $c$ to $d$ avoiding $\{a,b\}\cup A^\perp$. 
Let $c = x_1,\ldots, x_n =d$ be a shortest such path. 
Suppose $n = 2$.  
Then $\{a,b,c,d\}$ is a simplex. 
As $A$ is a maximal irreducible simplex, $d\in A^\perp$, 
which is a contradiction.  Hence $n > 2$. 

Now $b,x_1,\ldots,x_n$ is a cycle of length at least four. 
Hence the cycle $b,x_1,\ldots,x_n$ has a chord. 
The chord does not join $x_i$ to $x_j$ with $i < j$, 
since $x_1,\ldots, x_n$ is a shortest path from $c$ to $d$ avoiding $\{a,b\}\cup A^\perp$. 
Hence the chord joins $b$ to $x_i$ with $1 < i < n$. 
By repeating the same argument, starting with the cycles 
$b,x_1,\ldots, x_i$ and $b,x_i,\ldots,x_n$, 
we deduce that $m(b,x_i) < \infty$ for each $i$. 

Next consider the cycle $a,x_1,\ldots,x_n$. 
By the same argument as above, 
$m(a,x_i) < \infty$ for each $i$. 
Hence $\{a,b,c,x_2\}$ is a simplex. 
As $A$ is a maximal irreducible simplex, 
we must have that $x_2 \in A^\perp$, 
which is a contradiction. 
Thus $\{a,b\}\cup A^\perp$ is a $(c,d)$-separator. 
\end{proof}

Let $(W,S)$ be a Coxeter system. 
A set of {\it bad edges} for $(W,S)$ 
is a set ${\cal B}_2$ of pairs of elements $\{a,b\}$ of $S$ such that $m(a,b) \geq 5$ 
and if $A$ is an irreducible simplex of $(W,S)$ containing $\{a,b\}$ 
as a proper subset, then $A$ is of type ${\bf G}_3$ or ${\bf G}_4$. 
The elements of ${\cal B}_2$ are called {\it bad edges} of $(W,S)$. 
A {\it bad k-edge} is a bad edge $\{a,b\}$ with $m(a,b) = k$. 

Let ${\cal B}_2$ be a set of bad pairs for $(W,S)$. 
An irreducible simplex $A$ of $(W,S)$ is said to be {\it bad} 
if $A$ contains a bad pair. 
A bad irreducible simplex is of type ${\bf G}_3$, ${\bf G}_4$, 
or ${\bf D}_2(k)$ with $k\geq 5$. 
An irreducible simplex $A$ of $(W,S)$ is said to be {\it good} 
if $A$ is not bad. 
Let ${\cal B}$ be the set of all bad irreducible simplices of $(W,S)$. 
We call ${\cal B}$ a set of {\it bad irreducible simplices} for $(W,S)$. 

Let $(W,S)$ be a Coxeter system with a set ${\cal B}$ of bad irreducible simplices. 
A subset $B$ of $S$ is said to be a {\it bad separator} of $S$ if 
there is there is a pair of elements $\{a,b\}$ of $B$ with $m(a,b) = 2$, 
called a {\it pair of eyes} of $B$, 
and an element $c$ of $S-B$, called a {\it bad focus} of $B$, such that $A = \{a,b,c\}$, 
called a ${\it head}$ of $B$, is a bad maximal irreducible simplex of $(W,S)$ and 
$B \subset \{a,b\}\cup A^\perp$, 
and there is an element $f$ of $S-B$, called a {\it good focus} of $B$,  
such that $B$ is a minimal $(c,f)$-separator of $S$.

Let $B$ be a bad separator of $S$. 
If $\{a,b\}$ is a pair of eyes of $B$, 
then $\{a\}$ and $\{b\}$ are singleton components of $B$. 
If $A = \{a,b,c\}$ is a head of $B$, 
with eyes $\{a,b\}$, then $A$ is a component of $B\cup\{c\}$. 
If $f$ is a good focus of $B$, then $f \in S-(A\cup A^\perp)$. 

Consider the chordal Coxeter system $(W,S)$ with the P-diagram illustrated in Figure 1. 
In Figure 1, the edges incident with $e$ and $f$ can have any finite label. 
The set $\{a,b,c\}$ is a maximal irreducible simplex of $(W,S)$ of type ${\bf G}_3$.  
Hence we can declare the 5-edge $\{b,c\}$ to be bad. 
Then the set $\{a,b\}$ is a bad separator of $S$ with eyes $\{a,b\}$, 
bad focus $c$, and good focus $f$.   
Moreover, the set $\{a,b,d\}$ is a bad separator of $S$ with eyes $\{a,b\}$, 
bad focus $c$, and good focus $e$.

\medskip

$$\mbox{
\setlength{\unitlength}{.8cm}
\begin{picture}(8,8)(0,0)
\thicklines
\put(1,3){\circle*{.15}}
\put(1,3){\line(1,0){2.43}}
\put(7,3){\line(-1,0){3.1}}
\put(5,3){\circle*{.15}}
\put(7,3){\circle*{.15}}
\put(1,3){\line(1,1){4}}
\put(5,3){\line(0,1){4}}
\put(5,7){\circle*{.15}}
\put(1,3){\line(1,-1){2}}
\put(3,1){\line(1,1){3.16}}
\put(3,1){\line(1,3){2}}
\put(7,5){\line(-1,-1){.49}}
\put(7,5){\circle*{.15}}
\put(3,1){\line(1,3){2}}
\put(3,1){\line(2,1){4}}
\put(3,1){\circle*{.15}}
\put(5,7){\line(1,-1){2}}
\put(5,7){\line(1,-2){2}}
\put(3,3.3){5}
\put(2.6,5.25){3}
\put(5.2,4.7){2}
\put(4.9,7.4){$a$}
\put(4.9,2.35){$b$}
\put(.4,2.9){$c$}
\put(2.75,.25){$d$}
\put(7.4,2.9){$e$}
\put(7.4,4.9){$f$}
\put(1.5,1.5){2}
\put(4.2,3.6){2}
\put(3.8,2.2){2}
\end{picture}}$$

\medskip

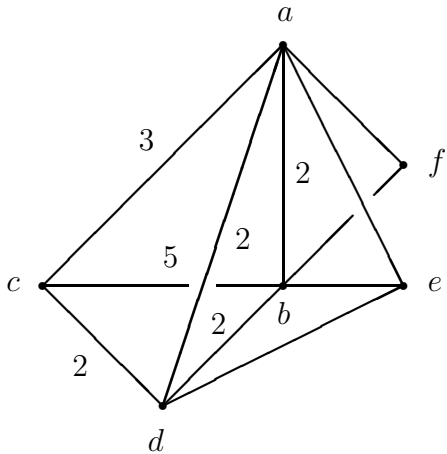
\centerline{\bf Figure 1}

\pagebreak

\begin{lemma} 
{\rm (Kloks and Kratsch \cite{K-K}, Lemma 2)}
Let $\Gamma$ be the P-diagram of a Coxeter system $(W,S)$, 
Let $B$ be a $(c,f)$-separator of $S$,   
and let $K_c$ and $K_f$ be the connected components of $\Gamma-B$ 
containing $c$ and $f$, respectively. 
Then $B$ is a minimal $(c,f)$-separator of $S$ if and only if each element of $B$ 
has a neighbor in both $K_c$ and $K_f$. 
\end{lemma}

\begin{lemma} 
Let $(W,S)$ be a chordal Coxeter system with a set ${\cal B}$ of bad irreducible simplices, 
and let $B$ be a subset of $S$,   
Then $B$ is a bad separator of $S$ if and only if there is 
a bad maximal irreducible simplex $A = \{a,b,c\}$ of $(W,S)$ with $m(a,b)=2$ 
such that $B \subset \{a,b\}\cup A^\perp$, 
and there is an element $d$ of $S-B$ such that $\{a,b,d\}$ is a simplex of $(W,S)$ 
and $B$ is a minimal $(c,d)$-separator of $S$.  
\end{lemma}
\begin{proof}
Suppose $B$ is a bad separator of $S$. 
Let $A = \{a,b,c\}$ be a head of $B$ with eyes $\{a,b\}$, 
and let $f$ be an element of $S-B$ such that $B$ is a minimal $(c,f)$-separator of $S$. 
By lemma 4.2, both $a$ and $b$ have neighbors $u$ and $v$, respectively, in the connected 
component $K_f$ of $\Gamma - B$ containing $f$. 
Let $u=x_1,\ldots,x_m=f$ and $f = y_1,\ldots, y_n=v$ be paths in $K_f$. 
Then $a,x_1,\ldots,x_m,y_2,\ldots,y_n,b$ is a cycle. 
As $\Gamma$ is chordal, there is a $d$ on this path  
such that $\{a,b,d\}$ is a simplex. 
As $d \in K_f$, we have that $K_f$ is the connected component $K_d$ of $\Gamma-B$ containing $d$.  
Therefore $B$ is a minimal $(c,d)$-separator of $S$ by Lemma 4.2. 

Conversely, suppose there is a bad maximal irreducible simplex $A = \{a,b,c\}$ 
of $(W,S)$ with $m(a,b)=2$ 
such that $B \subset \{a,b\}\cup A^\perp$, 
and there is an element $d$ of $S-B$ such that $\{a,b,d\}$ is a simplex of $(W,S)$ 
and $B$ is a minimal $(c,d)$-separator of $S$.  
As $\{a,b,c\}$ and $\{a,b,d\}$ are simplices, $\{a,b\}\subset B$, 
and so $B$ is a bad separator of $S$. 
\end{proof}

\begin{lemma} 
Let $(W,S)$ be a chordal Coxeter system with a set ${\cal B}$ of bad irreducible simplices.  
Then $S$ has a bad separator if and only if there is 
a bad maximal irreducible simplex $A = \{a,b,c\}$ of $(W,S)$ with $m(a,b)=2$ 
and there is a $d\in S-(A\cup A^\perp)$ such that $\{a,b,d\}$ is a simplex of $(W,S)$. 
\end{lemma}
\begin{proof}
Suppose $A = \{a,b,c\}$ is a bad maximal irreducible simplex of $(W,S)$ with $m(a,b) = 2$ 
and $d\in S-(A\cup A^\perp)$ such that $\{a,b,d\}$ is a simplex of $(W,S)$. 
By Lemma 4.1, we have that $\{a,b\}\cup A^\perp$ is a $(c,d)$-separator. 
Let $B \subset \{a,b\}\cup A^\perp$ be a minimal $(c,d)$-separator. 
Then $B$ is a bad separator of $S$ by Lemma 4.3. 
\end{proof}

\begin{lemma} 
Let $(W,S)$ be a chordal Coxeter system  
with a set ${\cal B}$ of bad irreducible simplices, 
and suppose that $S$ has no bad separators. 
Let $B$ be a simplex of $(W,S)$ 
with two singleton components $\{a\}$ and $\{b\}$ 
for which there is a $c\in S-B$ such that 
$A = \{a,b,c\}$ is a bad maximal irreducible simplex of $(W,S)$. 
Then $B-\{a,b\}\subset A^\perp$ and $B\cup\{c\}$ is a simplex of $(W,S)$. 
\end{lemma}
\begin{proof}
Let $d$ be an element of $B-\{a,b\}$. 
Then $\{a,b,d\}$ is a simplex. 
Hence $d\in A^\perp$ otherwise $S$ 
would have a bad separator by Lemma 4.4. 
Therefore $B-\{a,b\} \subset A^\perp$. 
As $B$ is a simplex, $B-\{a,b\}$ is a simplex, 
and so $B\cup\{c\} = A\cup (B-\{a,b\})$ is a simplex. 
\end{proof}

\begin{lemma}  
Let $(W,S)$ be a Coxeter system,  
and let $(S_1,S_0,S_2)$ be a separation of $S$ such that $S_0$ is a simplex.  
If $B$ is a minimal $(c,d)$-separator of $S_1$, 
then $B$ is a minimal $(c,d)$-separator of $S$. 
\end{lemma}
\begin{proof}
On the contrary, suppose $B$ is not a $(c,d)$-separator of $S$. 
Then there is a path $c=x_1,\ldots,x_n=d$ with distinct vertices 
in $\Gamma(W,S)$ that avoids $B$. 
If $x_i\in S_1$ for each $i$, we have a contradiction, 
since $B$ is a $(c,d)$-separator of $S_1$. 
Let $i$ be the first index such that $x_i\in S_2-S_0$, and let $j$ 
be the last index such that $x_j\in S_2-S_0$. 
Then $1 < i\leq j < n$ and $x_{i-1},x_{j+1} \in S_0$. 
As $S_0$ is a simplex, $c=x_1,\ldots,x_{i-1},x_{j+1},\ldots,x_n=d$ 
is a path in $\Gamma(W_1,S_1)$ that avoids $B$, which is a contradiction. 
Therefore $B$ is a $(c,d)$-separator of $S$. 

If we eliminate an element $e$ from $B$, 
then there is a path $c=x_1,\ldots,x_n=d$ in $\Gamma(W_1,S_1)$ 
that avoids $B-\{e\}$, since $B$ is a minimal $(c,d)$-separator of $S_1$. 
Hence $B$ is a minimal $(c,d)$-separator of $S$. 
\end{proof}

If $u$ is an element of a group $W$, 
we denote the inner automorphism $w \mapsto uwu^{-1}$ by $u_\ast$.  

\begin{theorem} 
{\rm (The Automorphic Angle Deformation Theorem)}
Let $(W,S)$ be a chordal Coxeter system of finite rank  
with a set ${\cal B}$ of bad irreducible simplices 
containing at most one bad 5-edge, and suppose that $S$ has no bad separators. 
For each $A$ in ${\cal B}$, let $\beta_A:\langle A\rangle \to \langle A\rangle$ 
be an automorphism that is not inner and that maps each element of $A$ 
to a conjugate of itself in $\langle A\rangle$. 
Then there is an automorphism $\alpha$ of $W$ 
such that for each simplex $B$ of $(W,S)$, 
there is an element $w_B$ of $W$ such that 

\pagebreak
\begin{enumerate}
\item $(w_B)_\ast\alpha\langle B\rangle = \langle B\rangle$, 
\item if $A$ is a bad component of $B$, 
then $(w_B)_\ast\alpha$ restricts on $\langle A\rangle$ to $\beta_A$,
\item if $A$ is a good component of $B$ of rank at least two, 
then $(w_B)_\ast\alpha$ restricts on $\langle A\rangle$ to the identity map, and  
\item if $A = \{a\}$ is a singleton component of $B$,  
then $(w_B)_\ast\alpha$ fixes $a$ unless 
there is another singleton component $\{b\}$ of $B$ 
and an element $c$ of $S-B$ such that $\{a,b,c\}$ is 
a bad maximal irreducible simplex of $(W,S)$ 
in which case $(w_B)_\ast\alpha$ transposes $a$ and $b$. 
\end{enumerate}
\end{theorem}
\begin{proof}
The proof is by induction on $|S|$. 
The theorem is clear if $|S|=1$, 
so assume $|S| > 1$ and the theorem is true for all chordal Coxeter systems 
of rank less than $|S|$. 
Assume first that $(W,S)$ is complete. 
Then $(W,S')$ is complete by Prop. 5.10 of \cite{M-R-T}. 
Let $(W,S) = (W_1,S_1)\times\cdots\times(W_n,S_n)$ be the factorization 
of $(W,S)$ into irreducible factors. 
Let $\alpha$ be the automorphism of $W$ 
that restricts on $W_i$ to the identity map if $S_i$ is good  
or to $\beta_i=\beta_{S_i}$ if $S_i$ is bad. 

Let $B$ be a simplex of $(W,S)$, 
and let $A$ be a component of $B$. 
Then $A \subset S_i$ for some index $i$. 
If $S_i$ is good, then $A$ is good, in which case 
define $w_i = 1$. 
Suppose $S_i$ is bad. 
Then $S_i$ is of type ${\bf G}_3$, ${\bf G}_4$ 
or ${\bf D}_2(k)$ with $k\geq 5$. 
If $A = S_i$, define $w_i = 1$. 
Now assume that $A$ is a proper subset of $S_i$. 
Then $A$ is of type ${\bf A}_1$, ${\bf A}_2$, ${\bf A}_3$, ${\bf D}_2(5)$, or ${\bf G}_3$. 
Let $B_i = B\cap S_i$. 
Then $B_i$ has one or two components. 
If $B_i$ has two components, then $A$ can be either component 
and $B_i$ is of type ${\bf A}_1\times{\bf A}_1$, ${\bf A}_1\times{\bf A}_2$, 
or ${\bf A}_1\times{\bf D}_2(5)$. 
By Prop. 32 of \cite{F-H}, 
there is an element $u_i$ of $W_i$ such that 
$(u_i)_\ast\beta_i\langle B_i\rangle$ is a visual subgroup of $(W_i,S_i)$. 
Moreover there is an element $v_i$ of $W_i$ such that 
$(v_i)_\ast(u_i)_\ast\beta_i\langle B_i\rangle = \langle B_i\rangle$ 
and $(v_iu_i)_\ast\beta_i$ leaves invariant each component of $B_i$ 
unless $S_i$ is of type ${\bf G}_3$ and $B_i = \{a,b\}$ 
with $m(a,b) = 2$, in which case $(v_iu_i)_\ast\beta_i$ transposes $a$ and $b$. 
If $A$ is of type ${\bf A}_2$ or ${\bf A}_3$, 
we may assume, by further conjugating by an element of $\langle A\rangle$,  
that $(v_iu_i)_\ast\beta_i$ restricts on $\langle A\rangle$ to the identity map. 
If $A$ is of type ${\bf D}_2(5)$ or ${\bf G}_3$, 
we may assume, by further conjugating by an element of $\langle A\rangle$,  
that $(v_iu_i)_\ast\beta_i$ restricts on $\langle A\rangle$ to $\beta_A$. 
Let $w_i = v_iu_i$. 
Then $w_B = w_1\cdots w_n$ has the desired properties 1-4. 

Now assume that $(W,S)$ is incomplete. 
Let $C \subset S$ be a minimal separator of $S$. 
Then $C$ is a simplex by Theorem 1 of \cite{Dirac}. 
Let $(S_1,C,S_2)$ be a separation of $S$, 
and let $W_i = \langle S_i\rangle$ for $i=1,2$. 
Then $(W_i,S_i)$ is a chordal Coxeter system for each $i=1,2$. 

Suppose $B$ is a bad separator of $S_1$. 
By Lemma 4.3, there is a bad maximal irreducible simplex $A = \{a,b,c\}$ of $(W_1,S_1)$ 
with $m(a,b) = 2$ such that $B\subset \{a,b\}\cup A^\perp$,  
and there is a $d \in S_1-B$ such that $\{a,b,d\}$ 
is a simplex of $(W_1,S_1)$ and $B$ is a minimal $(c,d)$-separator of $S_1$. 
Now $B$ is a minimal $(c,d)$-separator of $S$ by Lemma 4.6. 
As $S$ has no bad separators, $A$ is not a maximal irreducible simplex of $(W,S)$. 
Therefore there is a $e\in S_2-S_0$ such that $A\cup\{e\}$ is an irreducible simplex of $(W,S)$. 
Now $A\cup\{e\} \subset S_2$, since $A\cup\{e\}$ is a simplex. 
Hence $A \subset C$, and therefore $C$ contains the bad 5-edge, 
in which case we replace $C$ by a subset of $B$ that is a minimal separator of $S$. 
Then $C$ does not contain the bad 5-edge. 
Therefore $S_1$ and $S_2$ have no bad separators. 

By the induction hypothesis, 
there is an automorphism $\alpha_i$ of $W_i$ for each $i=1,2$ 
such that for each simplex $B$ of $(W_i,S_i)$, 
there is an element $w^{(i)}_B$ of $W_i$ such that 
(1) $(w^{(i)}_B)_\ast\alpha_i\langle B\rangle = \langle B\rangle$, 
(2) if $A$ is a bad component of $B$, 
then $(w^{(i)}_B)_\ast\alpha_i$ restricts on $\langle A\rangle$ to $\beta_A$,  
(3) if $A$ is a good component of $B$ of rank at least two, 
then $(w^{(i)}_B)_\ast\alpha_i$ restricts on $\langle A\rangle$ to the identity map, and  
(4) if $A = \{a\}$ is a singleton component of $B$, 
then $(w^{(i)}_B)_\ast\alpha_i$ fixes $a$ unless 
there is another singleton component $\{b\}$ of $B$ 
and an element $c$ of $S_i-B$ such that $\{a,b,c\}$ is 
a bad maximal irreducible simplex 
in which case $(w^{(i)}_B)_\ast\alpha_i$ transposes $a$ and $b$. 

Suppose $C$ has two singleton components $\{a\}$ and $\{b\}$ 
for which there is a $c\in S-C$ such that $A=\{a,b,c\}$ is a bad 
maximal irreducible simplex of $(W,S)$. 
Then $C\subset \{a,b\}\cup A^\perp$ by Lemma 4.5. 
We may assume $c\in S_1-C$. 
Let $f$ be an element of $S_2-C$. 
Then $C$ is a minimal $(c,f)$-separator of $S$. 
Therefore $C$ is a bad separator of $S$, 
which is a contradiction. 
Therefore $C$ does not have two singleton components $\{a\}$ and $\{b\}$ 
for which there is a $c\in S-C$ such that $A=\{a,b,c\}$ is a bad 
maximal irreducible simplex of $(W,S)$.  
Hence $(w^{(1)}_C)_\ast\alpha_1$ agrees with $(w^{(2)}_C)_\ast\alpha_2$ on $\langle C\rangle$. 
Let $\alpha$ be the automorphism of $W$ that restricts on $W_i$ to $(w^{(i)}_C)_\ast\alpha_i$ 
for each $i=1,2$. 

Let $B$ be a simplex of $(W,S)$. 
Then $B\subset S_i$ for some $i$. 
If $\{a\}$ and $\{b\}$ are singleton components of $B$ 
for which there is a $c\in S-B$ such that $\{a,b,c\}$ is 
a bad maximal irreducible simplex of $(W,S)$, 
then $B\cup \{c\}$ is a simplex by Lemma 4.5  
and we choose $i$ so that $B\cup \{c\} \subset S_i$. 
Let $w_B = w^{(i)}_B\big(w^{(i)}_C\big)^{-1}$. 
Then $w_B$ has the desired properties 1-4, 
which completes the induction. 
\end{proof}

\section{The Star Decomposition Theorem}

Let $(W,S)$ be a Coxeter system with a set ${\cal B}$ 
of bad irreducible simplices. 
A subset $D$ of $S$ is said to be a {\it gross separator} of $S$ if 
there is there is a pair of elements $\{a,b\}$ of $D$ with $m(a,b) = 2$, 
called a {\it pair of eyes} of $D$, 
and an element $c$ of $S-D$, called a {\it bad focus} of $D$, such that $A = \{a,b,c\}$, 
called a ${\it head}$ of $D$, is a bad maximal irreducible simplex of $(W,S)$ and 
$D = \{a,b\}\cup A^\perp$, 
and there is an element $f$ of $S-B$, called a {\it good focus} of $D$,  
such that $D$ is a $(c,f)$-separator of $S$
and neither $D-\{a\}$ nor $D-\{b\}$ is a $(c,f)$-separator of $S$. 

Let $D$ be a gross separator of $S$. 
If $\{a,b\}$ is a pair of eyes of $D$, 
then $\{a\}$ and $\{b\}$ are singleton components of $D$. 
If $A = \{a,b,c\}$ is a head of $D$, 
with eyes $\{a,b\}$, then $A$ is a component of $D\cup\{c\}$. 
If $f$ is a good focus of $D$, then $f \in S-(A\cup A^\perp)$. 

Consider the Coxeter system $(W,S)$ with the P-diagram illustrated 
in Figure 1 with $\{b,c\}$ a bad 5-edge. 
Observe that the set $\{a,b,d\}$ is a gross separator with eyes $\{a,b\}$, 
bad focus $c$, and good foci $e$ and $f$.

\begin{lemma} 
Let $(W,S)$ be a Coxeter system with just one bad 5-edge $\{x,y\}$. 
Let $C$ be a subset of $S$, and  
let $A = \{a,b,c\}$ be a bad irreducible simplex of $(W,S)$ with $m(a,b)=2$ 
such that $\{a,b\} \subset C\subset \{a,b\}\cup A^\perp$. 
Then the pair $\{a,b\}$ is unique and $\{c\} = \{x,y\}-\{a,b\}$.  
Moreover, if $C=\{a,b\}\cup A^\perp$, then $\{a,b\}$ determines $C$. 
\end{lemma}
\begin{proof}
Now $\{x,y\} \subset A$, since $A$ is bad. 
Hence $C-\{a,b\} \subset \{x,y\}^\perp$. 
Therefore $\{a,b\} = C-\{x,y\}^\perp$, and so the pair $\{a,b\}$ is unique. 
As $m(a,b) =2$, we have that $c=\{x,y\}-\{a,b\}$.  
If $C = \{a,b\}\cup A^\perp$, then $\{a,b\}$ determines $C$, 
since $\{a,b\}$ determines $A$. 
\end{proof} 

\begin{lemma} 
Let $(W,S)$ be a Coxeter system with just one bad 5-edge, 
and let $D$ be a gross separator of $S$ with head $A=\{a,b,c\}$ 
and eyes $\{a,b\}$. 
Then the set $F$ of good foci of $D$ 
is the set of vertices of all 
the connected components of $\Gamma(W,S)-D$ 
that contain a neighbor of $a$ and a neighbor of $b$ 
but do not contain $c$. 
\end{lemma}
\begin{proof}
Let $f$ be a good focus of $D$, and  
let $K_f$ be the connected component of $\Gamma-D$ 
that contains $f$.  Then $c \not\in K_f$, 
since $D$ is a $(c,f)$-separator.  

Now $D-\{a\}$ is not a $(c,f)$-separator, and so 
there is a path $c=x_1,\ldots,x_n =f$ in $\Gamma$ with distinct vertices that avoids $D-\{a\}$. 
The path cannot avoid $D$, since $D$ is a $(c,f)$-separator, 
and so $a = x_i$ for some $i$ with $1 < i < n$. 
Then $x_{i+1}$ is a neighbor of $a$ in $K_f$. 
Likewise $K_f$ contains a neighbor of $b$. 

Let $g$ be an element of $S-D$ and let $K_g$ 
be the connected component of $\Gamma-D$ that contains $g$. 
Suppose that $K_g$ contains a neighbor $u$ of $a$ and a neighbor 
of $v$ of $b$ but does not contain $c$. 
Then $D$ is a $(c,g)$-separator of $S$. 
Let $u=x_1,\ldots,x_n=g$ be a path in $K_g$. 
Then $c,a,x_1,\ldots,x_n$ is a path from $c$ to $g$ avoiding $D-\{a\}$,  
and so $D-\{a\}$ is not a $(c,g)$-separator of $S$. 
Likewise $D-\{b\}$ is not a $(c,g)$-separator of $S$. 
Therefore $g$ is a good focus of $D$. 
Thus $F$ is the set of vertices of all 
the connected components of $\Gamma-D$ 
that contain a neighbor of $a$ and a neighbor of $b$ 
but do not contain $c$. 
\end{proof}

\begin{lemma} 
Let $(W,S)$ be a Coxeter system with  
just one bad 5-edge $\{x,y\}$. 
Let $(S_1,D_1,S_2)$ be a separation of $S$ 
with $D_1$ a gross separator of $S$ 
and $S_1-D_1$ the set of all good foci of $D_1$. 
If $D_2$ is a gross separator of $S_2$, 
then $D_2$ is a gross separator of $S$ 
and $f\in S-D_2$ is a good focus of $D_2$ relative to $S$ 
if and only if $f \in S_2-D_2$ 
and $f$ is a good focus of $D_2$ relative to $S_2$. 
\end{lemma}
\begin{proof}
Let $A_i=\{a_i,b_i,c_i\}$ be the head of $D_i$ with eyes $\{a_i,b_i\}$ 
for each $i=1,2$. 
We claim that $D_1\neq D_2$. 
On the contrary, assume that $D_1=D_2$. 
Then $A_1=A_2$ by Lemma 5.1. 
Let $f_2 \in S_2-D_2$ be a good focus of $D_2$ relative to $S_2$. 
Then neither $D_2-\{a_2\}$ nor $D_2-\{b_2\}$ is a $(c_2,f_2)$-separator of $S_2$. 
Hence neither $D_2-\{a_2\}$ nor $D_2-\{b_2\}$ is a $(c_2,f_2)$-separator of $S$. 
Now $f_2$ is not a good focus of $D_2$ relative to $S$, since $f_2 \in S_2-D_1$. 
Hence $D_2$ is not a $(c_2,f_2)$-separator of $S$, and so 
there is a path $c_2=x_1,\ldots,x_n=f_2$ in $\Gamma(W,S)$ avoiding $D_2$. 
The path $x_1,\ldots,x_n$ must lie in $S_2-D_2$, since $(S_1,D_2,S_2)$ is a separation of $S$. 
Therefore $D_2$ is not a $(c_2,f_2)$-separator of $S_2$, which is a contradiction. 
Thus $D_1\neq D_2$.

Let $e$ be an element of $S-A_2$ 
such that $A_2\cup\{e\}$ is a simplex of $(W,S)$. 
We claim that $e\in S_2$. 
On the contrary, suppose that $e \in S-S_2$. 
Then $e\in S_1-D_1$.  
Hence $A_2\cup\{e\} \subset S_1$, 
and so $A_2\subset D_1$,  
whence $\{x,y\} \subset D_1$, which is a contradiction. 
Therefore $e \in S_2$. 
Hence $A_2$ is a maximal irreducible simplex of $(W,S)$ 
and $D_2 = \{a_2,b_2\}\cup A_2^\perp$ in $(W,S)$. 

Let $f_2 \in S_2-D_2$ be a good focus of $D_2$ relative to $S_2$. 
We claim that $D_2$ is a $(c_2,f_2)$-separator of $S$. 
On the contrary, suppose $D_2$ is not a $(c_2,f_2)$-separator of $S$. 
Then there is a path $c_2=x_1,\ldots, x_n = f_2$ with distinct vertices 
in $\Gamma(W,S)$ that avoids $D_2$. 
If $ x_i\in S_2$ for each $i$, we have a contradiction, 
since $D_2$ is a $(c_2,f_2)$-separator of $S_2$. 
Let $i$ be the first index such that $x_i\in S_1-D_1$ and let $j$ 
be the last index such that $x_j\in S_1-D_1$. 
Then $1< i\leq j< n$ and $x_{i-1},x_{j+1}\in D_1$. 
If $x_{i-1}$ and $x_{j+1}$ are adjacent in $\Gamma(W,S)$, 
then $c_2 = x_1,\ldots,x_{i-1},x_{j+1},\ldots, x_n = f_2$ 
is a path in $S_2$ that avoids $D_2$, which is a contradiction. 
Therefore $x_{i-1},x_{j+1} \in A_1^\perp-\{a_1,b_1\}$. 
As $D_1\neq D_2$, we have that $\{a_1,b_1\}\neq\{a_2,b_2\}$ by Lemma 5.1. 
Say $a_1 \not\in\{a_2,b_2\}$. 
Now $a_1 \not\in A_2^\perp$, since $A_1$ and $A_2$ share the common 5-edge $\{x,y\}$  
and $A_1$ is irreducible. 
Therefore $a_1 \not\in D_2$. 
Hence $c_2=x_1,\ldots,x_{i-1},a_1,x_{j+1},\ldots,x_n=f_2$ 
is a path in $S_2$ that avoids $D_2$, which is a contradiction. 
Therefore $D_2$ is a $(c_2,f_2)$-separator of $S$. 
Moreover neither $D_2-\{a_2\}$ nor $D_2-\{b_2\}$ 
is a $(c_2,f_2)$-separator of $S$, 
since neither $D_2-\{a_2\}$ nor $D_2-\{b_2\}$ 
is a $(c_2,f_2)$-separator of $S_2$. 
Thus $D_2$ is a gross separator of $S$.  

Let $f\in S-D_2$ be a good focus of $D_2$ relative to $S$. 
We claim that $f\in S_2-D_2$. 
On the contrary, suppose $f \in S-S_2$, 
Then $f\in S_1-D_1$. 
Hence $f$ is a good focus of $D_1$. 
By Lemma 5.2, there is a path 
$f=x_1,\ldots,x_n=a_1$ in $\Gamma(W,S)$ with $x_1,\ldots,x_{n-1}$ in $S_1-D_1$. 
We have that $a_1\not\in D_2$ and $D_2 \subset S_2$. 
Hence $x_1,\ldots, x_n$ is a path from $f$ to $a_1$ that avoids $D_2$. 
If $a_1 = c_2$, then $x_1,\ldots, x_n$ is a path from $f$ to $c_2$ that avoids $D_2$, 
which is a contradiction, since $f$ is a focus of $D_2$. 
If $a_1\neq c_2$, then $a_1$ and $c_2$ are adjacent in $\Gamma(W,S)$, 
since $A_1$ and $A_2$ share the common 5-edge $\{x,y\}$. 
Hence $x_1,\ldots,x_n,c_2$ is a path from $f$ to $c_2$ that avoids $D_2$, 
which is contradiction. 
Therefore $f \in S_2-D_2$. 

As $D_2$ is a $(c_2,f)$-separator of $S$, 
we have that $D_2$ is a $(c_2,f)$-separator of $S_2$. 
Now $D_2-\{a_2\}$ is not a $(c_2,f)$-separator of $S$. 
Let $c_2=x_1,\ldots,x_n=f$ be a path in $\Gamma(W,S)$ that avoids $D_2-\{a_2\}$. 
If $x_i\in S_2$, for each $i$, then $x_1,\ldots,x_n$ is a 
path in $S_2$ that avoids $D_2-\{a_2\}$. 
Suppose the path $x_1,\ldots, x_n$ does not lie in $S_2$. 
Let $i$ be the first index such that $x_i\in S_1-D_1$ 
and let $j$ be the last index such that $x_j\in S_1-D_1$. 
Then $1< i\leq j< n$ and $x_{i-1},x_{j+1}\in D_1$. 
If $x_{i-1}$ and $x_{j+1}$ are adjacent in $\Gamma(W,S)$, 
then $c_2 = x_1,\ldots,x_{i-1},x_{j+1},\ldots, x_n = f$ 
is a path in $S_2$ that avoids $D_2-\{a_2\}$. 
If $x_{i-1}$ and $x_{j+1}$ are not adjacent in $\Gamma(W,S)$, 
then $c_2 = x_1,\ldots,x_{i-1},a_1,x_{j+1},\ldots, x_n = f$ 
is a path in $S_2$ that avoids $D_2-\{a_2\}$. 
Therefore $D_2-\{a_2\}$ is not a $(c_2,f)$-separator of $S_2$. 
Likewise $D_2-\{b_2\}$ is not a $(c_2,f)$-separator of $S_2$. 
Therefore $f$ is a good focus of $D_2$ relative to $S_2$. 
\end{proof}

\pagebreak
\begin{lemma} 
Let $(W,S)$ be a Coxeter system of finite rank 
with just one bad 5-edge $\{x,y\}$.  
Then $(W,S)$ has a unique reduced visual graph of groups decomposition $\Upsilon$ 
whose underlying graph is a star such that if $(V_0,S_0)$ is the center 
vertex system of $\Upsilon$, then $\{x,y\}\subset S_0$ 
and $S_0$ has no gross separators, 
and if $(E_1,D_1),\ldots,(E_m,D_m)$ are the edge systems of $\Upsilon$, 
then $D_1,\ldots,D_m$ are the gross separators of $(W,S)$, 
and if $(V_1,S_1),\ldots, (V_m,S_m)$ are the vertex systems of $\Upsilon$ 
such that $D_i = S_0\cap S_i$ for each $i=1,\ldots,m$, 
then $S_i-D_i$ is the set of all the good foci of $D_i$ for each $i=1,\ldots,m$. 
\end{lemma}
\begin{proof}
The proof is by induction on $|S|$. 
The theorem is clear if $S =\{x,y\}$, 
so assume $|S|>2$ and the theorem is true 
for all Coxeter systems of rank less than $|S|$. 
If $S$ has no gross separators, then the trivial graph of groups decomposition $\Upsilon$ of $(W,S)$ 
with only one vertex group has the desired properties. 

Suppose $S$ has a gross separator $D_1$. 
Let $F_1$ be the set of all good foci of $D_1$, 
let $S_1 = D_1\cup F_1$, and let $T_1 = D_1\cup (S-F_1)$.   
Then $(S_1,D_1,T_1)$ is a separation of $S$ such that $\{x,y\}\subset T_1$ by Lemma 5.2.  
By the induction hypothesis, 
the system $(\langle T_1\rangle,T_1)$ has a reduced visual 
graph of groups decomposition $\Upsilon_1$ whose underlying graph 
is a star such that if $(V_0,S_0)$ is the center vertex system of $\Upsilon_1$, 
then $\{x,y\}\subset S_0$ and $S_0$ has no gross separators, 
and if $(E_2,D_2),\ldots,(E_{m},D_{m})$ are the edge systems of $\Upsilon_1$, 
then $D_2,\ldots,D_{m}$ are the gross separators of $T_1$, 
and if $(V_2,S_2),\ldots, (V_{m},S_{m})$ are the vertex systems of $\Upsilon_1$ 
such that $D_i=S_0\cap S_i$ for $i=2,\ldots,m$, then $S_i-D_i$ 
is the set $F_i$ of all good foci of $D_i$ relative to $T_1$ for each $i=2,\ldots,m$.

Let $A_i = \{a_i,b_i,c_i\}$ be the head of $D_i$ with eyes $\{a_i,b_i\}$ 
for each $i=1,\ldots,m$.  
Now if $e \in A_1^\perp$, then $A_1\cup\{e\}$ is a simplex, 
and so $A_1\cup\{e\} \subset S_0$, since $V_0$ is the only vertex group of $\Upsilon_1$ 
that contains the bad edge $\{x,y\}$. 
Therefore $D_1\subset S_0$. 
Hence we have a reduced visual graph of groups decomposition $\Upsilon$ of $(W,S)$ 
whose underlying graph is a star such that $(V_0,S_0)$ is the center 
vertex system, $(E_1,D_1), \ldots, (E_m,D_m)$ are the edge systems,
and $(V_1,S_1),\ldots,(V_m,S_m)$ are the vertex systems 
such that $D_i=S_0\cap S_i$ for each $i=1,\ldots,m$. 
We have that $\{x,y\}\subset S_0$ and $S_0$ has no gross separators. 
By Lemma 5.3, we have that $D_2,\ldots,D_m$ are gross separators of $S$ 
and $F_i=S_i-D_i$ is the set of all the good foci of $D_i$ relative to $S$ 
for each $i=2,\ldots,m$. Therefore $D_1 \neq D_i$ for all $i=2,\ldots,m$. 

Let $D$ be a gross separator of $S$, 
and let $A=\{a,b,c\}$ be the head of $D$ with eyes $\{a,b\}$. 
We claim that $D=D_i$ for some $i$. 
On the contrary, suppose $D\neq D_i$ for each $i$. 
Now $D\subset S_0$ by the same argument that we used to show that $D_1\subset S_0$. 
Let $f$ be a good focus for $D$. 
We claim that $f\in S-S_0$. 
On the contrary, suppose $f\in S_0$. 
As $D$ is a $(c,f)$-separator of $S$, 
we have that $D$ is a $(c,f)$-separator of $S_0$. 
Now $S_0$ has no gross separators, 
and so either $D-\{a\}$ or $D-\{b\}$ is a $(c,f)$-separator of $S_0$. 
Say $D-\{a\}$ is a $(c,f)$-separator of $S_0$.
Now $D-\{a\}$ is not a $(c,f)$-separator of $S$. 
Hence there is a shortest path $c=x_1,\ldots,x_n=f$ in $\Gamma(W,S)$ that avoids $D-\{a\}$. 
If $x_i\in S_0$ for each $i=1,\ldots,n$, we have a contradiction. 
Suppose $i$ is the first index such that $x_i\not\in S_0$.  
Then $x_i\in S_k-D_k$ for some $k > 0$.  
Let $j$ be the last index such that $x_j\in S_k-D_k$. 
Then $1< i\leq j< n$ and $x_{i-1},x_{j+1}\in D_k$. 
Now $x_{i-1}$ and $x_{j+1}$ are not adjacent in $\Gamma(W,S)$. 
As in the proof of Lemma 5.3, there is an eye of $D_k$, say $a_k$, that is not in $D$. 
Then $c=x_1,\ldots,x_{i-1},a_k,x_{j+1},\ldots,x_n = f$ is a shortest path. 
After applying all such replacements, we may assume that $x_i\in S_0$ for each $i$, 
which is a contradiction. 
Therefore $f \in S-S_0$. 

Now $f\in S_k-S_0$ for some $k$. 
Then $f$ is a focus of $D_k$. 
As before, $D_k$ has an eye, say $a_k$, that is not in $D$. 
By Lemma 5.2, there is a path $f=x_1,\ldots,x_n=a_k$ 
in $\Gamma(W,S)$ with $x_1,\ldots,x_{n-1}$ in $S_k-D_k$. 
Hence $x_1,\ldots, x_n$ is a path from $f$ to $a_k$ that avoids $D$. 
If $a_k = c$, then $x_1,\ldots, x_n$ is a path from $f$ to $c$ that avoids $D$, 
which is a contradiction, since $f$ is a focus of $D$. 
If $a_k\neq c$, then $a_k$ and $c$ are adjacent in $\Gamma(W,S)$, 
since $A_k$ and $A$ share the common 5-edge $\{x,y\}$. 
Hence $x_1,\ldots,x_n,c$ is a path from $f$ to $c$ that avoids $D$, 
which is contradiction. 
Therefore $D = D_i$ for some $i$. 
Hence $D_1,\ldots,D_m$ are the gross separators of $S$. 

Finally, $\Upsilon$ is unique, since the edge and noncentral vertex systems 
are determined by the set of gross separators of $S$, 
and the center vertex system $(V_0,S_0)$ is determined 
by the relation $S_0 = S-(F_1\cup \cdots \cup F_m)$. 
\end{proof}

Let $(W,S)$ be a Coxeter system, let $c$ be an element of $S$, 
and let $B$ be a subset of $S-\{c\}$.  
Then $B$ is said to be {\it close} to $c$ 
if every element of $B$ is adjacent to $c$ in $\Gamma(W,S)$. 

\begin{lemma} 
{\rm (Kloks and Kratsch \cite{K-K}, Lemma 4)} 
Let $(W,S)$ be a Coxeter system, 
and let $c$ and $f$ be nonadjacent vertices of $\Gamma(W,S)$. 
Then there exists exactly one minimal $(c,f)$-separator of $S$ close to $c$. 
\end{lemma}

\begin{lemma} 
Let $(W,S)$ be a Coxeter system with just one bad 5-edge $\{x,y\}$.    
Let $D$ be a gross separator of $S$ with head $A$, 
and let $f$ be a good focus of $D$. 
Then $D$ contains a bad separator $B$ of $S$ 
with head $A$ and good focus $f$. 
Moreover $B$ is the unique bad separator of $S$ 
with head $A$ and good focus $f$. 
\end{lemma}
\begin{proof}
Let $A = \{a,b,c\}$ with $\{a,b\}$ the eyes of $D$. 
Now $D$ is a $(c,f)$-separator of $S$.  
Therefore $D$ contains a minimal $(c,f)$-separator $B$ of $S$. 
As neither $D-\{a\}$ nor $D-\{b\}$ is a $(c,f)$-separator of $S$, 
we deduce that $\{a,b\} \subset B$. 
Hence $B$ is a bad separator of $S$ with head $A$ and good focus $f$. 

Let $B$ be a bad separator of $S$ with head $A$ and good focus $f$. 
Then $B$ is a minimal $(c,f)$-separator of $S$ that is close to $\{c\}$, 
since $B\subset \{a,b\}\cup A^\perp$. 
Therefore $B$ is the unique bad separator of $S$ with head $A$ and good focus $f$ by Lemma 5.5. 
\end{proof}

\begin{lemma} 
Let $(W,S)$ be a Coxeter system with just one bad 5-edge $\{x,y\}$.    
Let $B$ be a bad separator of $S$ with head $A = \{a,b,c\}$ and eyes $\{a,b\}$. 
Let $D = \{a,b\}\cup A^\perp$. 
Let $f \in S-B$ be a good focus of $B$, 
and let $K_f$ be the connected component of $\Gamma(W,S)-B$ that contains $f$. 
Then $K_f$ is the connected component $C_f$ of $\Gamma(W,S)-D$ that contains $f$. 
\end{lemma}
\begin{proof}
First of all, $f \in S-D$, since $c$ and $f$ are nonadjacent in $\Gamma(W,S)$. 
Now $C_f \subset \Gamma-D \subset \Gamma-B$, and so $C_f \subset K_f$. 
We claim that $C_f = K_f$. 
On the contrary, let $e$ be a vertex in $K_f$ with $e$ adjacent to a vertex $g \in C_f$. 
Then $e \in D$. 
As $\{a,b\} \subset B$ and $e\in \Gamma-B$, we have that $e\in A^\perp$. 
Hence $e$ is adjacent to $c$.  
But no vertex of $K_f$ is adjacent to $c$, 
since $c$ is in a different component of $\Gamma-B$. 
Therefore $C_f$ has the same vertices as $K_f$. 

Let $e \in S-K_f$ be adjacent to a vertex $g \in K_f$, 
Then $e \in B$, and so $e \in D$ 
and $e$ is adjacent to the vertex $g \in C_f$. 
Therefore $C_f = K_f$ as subsets of $\Gamma$. 
\end{proof}

\begin{lemma} 
Let $(W,S)$ be a Coxeter system with just one bad 5-edge $\{x,y\}$.    
Let $B$ be a bad separator of $S$ with head $A = \{a,b,c\}$, eyes $\{a,b\}$, 
and good focus $f$.  
Let $D = \{a,b\}\cup A^\perp$. 
Then $D$ is a gross separator of $S$ with head $A$ and good focus $f$. 
\end{lemma}
\begin{proof}
As $B$ is a $(c,f)$-separator of $S$, 
we have that $D$ is a $(c,f)$-separator of $S$. 
Now neither $B-\{a\}$ nor $B-\{b\}$ is a $(c,f)$-separator of $S$. 
Let $K_f$ be the component of $\Gamma-B$ that contains $f$. 
Then $K_f$ contains a neighbor of $a$ and a neighbor of $b$ by Lemma 4.2. 
Now $f \in S-D$. 
Let $C_f$ be the component of $\Gamma-D$ that contains $f$. 
Then $C_f = K_f$ by Lemma 5.7. 
Hence $C_f$ contains a neighbor of $a$ and a neighbor of $b$. 
Therefore neither $D-\{a\}$ nor $D-\{b\}$ is a $(c,f)$-separator of $S$. 
Thus $D$ is a gross separator of $S$ with head $A$ and good focus $f$.   
\end{proof}

\begin{lemma} 
Let $(W,S)$ be a Coxeter system with just one bad 5-edge $\{x,y\}$. 
Let $D$ be a gross separator of $S$ with head $A$, 
and let $F$ be the set of all the good foci of $D$. 
Let $B_1,\ldots, B_k$ be the bad separators of $S$ with head $A$, 
and let $F_i$ be the set of all the good foci of $B_i$ for each $i$. 
Then $F = F_1\cup \cdots \cup F_k$. 
The sets $F_1,\ldots, F_k$ are pairwise disjoint, 
and if $i\neq j$, then no element of $F_i$ is adjacent to an element of $F_j$ 
in $\Gamma(W,S)$. 
\end{lemma}
\begin{proof}
Let $f$ be an element of $F$. 
Then $f\in F_i$ for exactly one index $i$ by Lemma 5.6. 
Let $f$ be an element of $F_i$ for some index $i$. 
Then $f \in F$ by Lemma 5.8. 
Therefore $F$ is the disjoint union of the sets $F_1,\ldots,F_k$. 

Let $A=\{a,b,c\}$ with $\{a,b\}$ the pair of eyes of $D$. 
Now $F_i$ is the set of vertices of all the connected components 
$C$ of $\Gamma-D$ such that each element of $B_i$ has a neighbor in $C$ 
and $C$ does not contain $c$ by Lemmas 4.2 and 5.7 for each $i$. 
Therefore if $i\neq j$, then no element of $F_i$ is adjacent to an element of $F_j$ 
in $\Gamma$. 
\end{proof}

\begin{theorem} 
{\rm (The Star Decomposition Theorem)} 
Let $(W,S)$ be a Coxeter system of finite rank 
with just one bad 5-edge $\{x,y\}$.  
Then $(W,S)$ has a unique reduced visual graph of groups decomposition $\Psi$ 
whose underlying graph is a star such that if $(V_0,S_0)$ is the center 
vertex system of $\Psi$, then $\{x,y\}\subset S_0$ 
and $S_0$ has no bad separators, 
and if $(E_1,B_1),\ldots,(E_n,B_n)$ are the edge systems of $\Psi$, 
then $B_1,\ldots,B_n$ are the bad separators of $(W,S)$, 
and if $(V_1,S_1),\ldots, (V_n,S_n)$ are the vertex systems of $\Psi$ 
such that $B_i = S_0\cap S_i$ for each $i=1,\ldots,n$, 
then $S_i-B_i$ is the set $F_i$ of all the good foci of $B_i$ for each $i=1,\ldots,n$.  
\end{theorem}
\begin{proof}
Let $\Upsilon$ be the reduced visual graph of groups decomposition of $(W,S)$ 
in Lemma 5.4. 
Then $\Upsilon$ refines to a reduced visual graph of groups decomposition $\Psi$ of $(W,S)$ 
as in the statement of the theorem by Lemmas 5.6, 5.8, and 5.9. 
Moreover, the center vertex system $(V_0,S_0)$ of $\Psi$ 
is the center vertex system of $\Upsilon$. 

Finally, $\Psi$ is unique, since the edge and noncentral vertex systems 
are determined by the set of bad separators of $S$, 
and the center vertex system $(V_0,S_0)$ is determined 
by the relation $S_0 = S-(F_1\cup \cdots \cup F_n)$. 
\end{proof}

\section{The Angle Deformation Theorem}

\begin{lemma} 
Let $(W,S)$ be a Coxeter system, 
and let $(S_1,S_0,S_2)$ be a separation of $S$. 
If $S_i'$ is a set of Coxeter generators for $\langle S_i\rangle$ for each $i=0,1,2$ 
such that $S_i'$ contains $S_0'$ for each $i = 1,2$,  
then $S'=S_1'\cup S_2'$ is a set of Coxeter generators for $W$. 
\end{lemma}
\begin{proof}
This follows from the fact that $W = \langle S_1\rangle\ast_{\langle S_0\rangle}\langle S_2\rangle$ 
is an amalgamated product decomposition. 
\end{proof}

\begin{theorem} 
{\rm (The 5-Edge Angle Deformation Theorem)} 
Let $(W,S)$ be a chordal Coxeter system of finite rank 
with just one bad edge $\{x,y\}$ 
and suppose $m(x,y) = 5$. 
Let ${\cal B}$ be the set of bad irreducible simplices of $S$, 
and for each $A$ in ${\cal B}$, let $\beta_A:\langle A\rangle \to \langle A\rangle$ 
be a reflection preserving automorphism that is not inner. 
Suppose that $S$ has a bad separator.  
Let $B_1,\ldots, B_n$ be the bad separators of $S$, 
let $\{a_i,b_i\}$ be the pair of eyes of $B_i$ for each $i$, 
and let $F_i$ be the set of good foci of $B_i$ for each $i$.  

Then $W$ has a set of Coxeter generators $S'$ such that for each $A$ in ${\cal B}$, 
there is a $w_A\in W$ such that 
$w_A\beta_A(A)w_A^{-1} \subset S'$,  
and for each good irreducible simplex $A$ of $(W,S)$, 
there is a $w_A\in W$ such that $w_A A w_A^{-1} \subset S'$. 

Moreover $\Gamma(W,S')$ is isomorphic to the labeled diagram $\Gamma'(W,S)$ 
obtained from $\Gamma(W,S)$ 
by cutting each edge $\{a_i,d\}$ with $d\in F_i$ at the vertex $a_i$ and 
then pasting the vertex $a_i$ of $\{a_i,d\}$ to the vertex $b_i$, 
and cutting each edge $\{b_i,e\}$ with $e\in F_i$ at the vertex $b_i$, 
and then pasting the vertex $b_i$ of $\{b_i,e\}$ to the vertex $a_i$ 
for each $i=1,\ldots,n$. 
\end{theorem}
\begin{proof}
Let $\Psi$ be the reduced visual graph of groups decomposition of $(W,S)$ 
in the star decomposition theorem. 
By the automorphic angle deformation theorem, 
there is an automorphism $\alpha$ of $V_0$ such that 
for each simplex $B$ of $(V_0,S_0)$, there is an element $u_B$ of $V_0$ such that 
\begin{enumerate}
\item $(u_B)_\ast\alpha\langle B\rangle = \langle B\rangle$, 
\item if $A$ is a bad component of $B$, 
then $(u_B)_\ast\alpha$ restricts on $\langle A\rangle$ to $\beta_A$,
\item if $A$ is a good component of $B$ of rank at least two, 
then $(u_B)_\ast\alpha$ restricts on $\langle A\rangle$ to the identity map, and  
\item if $A = \{a\}$ is a singleton component of $B$, 
then $(u_B)_\ast\alpha$ fixes $a$ unless 
there is another singleton component $\{b\}$ of $B$ 
and an element $c$ of $S_0-B$ such that $\{a,b,c\}$ is 
a bad maximal irreducible simplex of $(W,S)$ 
in which case $(u_B)_\ast\alpha$ transposes $a$ and $b$. 
\end{enumerate}
Let $R_0 = \alpha(S_0)$. 
Then $R_0$ is a set of Coxeter generators of $V_0$. 

Now $B_i$ is a simplex for each $i$ by Theorem 1 of Dirac \cite{Dirac}. 
Let $u_i = u_{B_i}$ for each $i$. 
Then $(u_i)_\ast\alpha(B_i) = B_i$ for each $i$, 
since $B_i$ has no bad components. 
Hence $u_i^{-1}B_iu_i = \alpha(B_i) \subset R_0$ 
and $B_i \subset u_iR_0 u_i^{-1}$ for each $i$. 

Let $W_i = \langle S_0\cup S_1\cup \cdots\cup S_i\rangle$ for $i=1,\ldots,n$. 
Then $R_1 = u_1R_0u_1^{-1}\cup S_1$ is a set of Coxeter generators for $W_1$ by Lemma 6.1 
applied to the separation $(S_0,B_1,S_1)$ of $S_0\cup S_1$. 
Likewise 
$$R_2 = u_2u_1^{-1}R_1u_1u_2^{-1}\cup S_2 = u_2R_0u_2^{-1}\cup u_2u_1^{-1}S_1u_1u_2^{-1}\cup S_2$$ 
is a set of Coxeter generators for $W_2$ by Lemma 6.1 
applied to the separation $(S_0\cup S_1,B_2,S_2)$ of $S_0\cup S_1\cup S_2$. 
Continuing in this way, we deduce that 
$$R_n = u_nR_0u_n^{-1}\cup u_nu_1^{-1}S_1u_1u_n^{-1}\cup\cdots\cup u_nu_{n-1}^{-1}S_{n-1}
u_{n-1}u_n^{-1}\cup S_n$$ 
is a set of Coxeter generators for $W$. 
Define 
$$S'=u_n^{-1}R_n u_n = R_0\cup u_1^{-1}S_1u_1\cup\cdots\cup u_n^{-1}S_nu_n.$$
Then $S'$ is a set of Coxeter generators for $W$. 

Let $A$ be a bad irreducible simplex of $(W,S)$. 
Then $A\subset S_0$, since $\{x,y\} \subset S_0$ and $\{x,y\}\not\subset B_i$ for each $i$. 
Now $(u_A)_\ast\alpha(A) = \beta_A(A)$, 
and so $u_A^{-1}\beta_A(A)u_A=\alpha(A)$. 
Let $w_A = u_A^{-1}$. 
Then $w_A\beta_A(A)w_A^{-1} \subset S'$. 

Let $A$ be a good irreducible simplex of $(W,S)$. 
If $A \subset S_0$, then we have $(u_A)_\ast\alpha(A) = A$, 
and so $u_A^{-1}Au_A = \alpha(A)$. 
Let $w_A = u_A^{-1}$. 
Then $w_AAw_A^{-1} \subset S'$. 
Now suppose $A \not\subset S_0$. 
Then $A\subset S_i$ for some $i>0$. 
Hence $u_i^{-1}Au_i \subset u_i^{-1}S_iu_i$. 
Let $w_A = u_i^{-1}$. 
Then $w_AAw_A^{-1} \subset S'$. 

Define $\phi:\Gamma'(W,S) \to \Gamma(W,S')$ so that  if $s\in S_0$, 
then $\phi(s) = \alpha(s)$ and if $s\in S_i-S_0$, then $\phi(s) = u_i^{-1}su_i$. 
Let $\{a,b\}$ be an edge of $\Gamma(W,S)$ with $a\in S_0$ and $b\in S_i-S_0$. 
Then $a\in B_i$, since $(S_0\cup (S-S_i), B_i, S_i)$ is a separation of $S$. 
Assume first that $a \in \{a_i,b_i\}$. 
Let $\ov a_i = b_i$ and $\ov b_i=a_i$. 
Then $\{\ov a,b\}$ is the corresponding edge of $\Gamma'(W,S)$ 
with edge label $m(a,b)$. 
The edge $\{\phi(\ov a),\phi(b)\} = \{\alpha(\ov a), u_i^{-1}bu_i\}$ 
of $\Gamma(W,S')$ has the same edge label, since 
$$\{\alpha(\ov a),u_i^{-1}bu_i\}=u_i^{-1}\{u_i\alpha(\ov a)u_i^{-1},b\}u_i=u_i^{-1}\{a,b\}u_i.$$
Now assume $a\not\in\{a_i,b_i\}$. 
Then $\{a,b\}$ is an edge of $\Gamma'(W,S)$ with edge label $m(a,b)$. 
The edge $\{\phi(a),\phi(b)\} = \{\alpha(a),u_i^{-1}bu_i\}$ of $\Gamma(W,S')$ 
has the same edge label, since
$$\{\alpha(a),u_i^{-1}bu_i\}=u_i^{-1}\{u_i\alpha(a)u_i^{-1},b\}u_i=u_i^{-1}\{a,b\}u_i.$$
Therefore $\phi$ is an isomorphism of labeled diagrams. 
\end{proof}

We say that the diagram $\Gamma'(W,S)$ in Theorem 6.2 
has been obtained from the diagram $\Gamma(W,S)$ by the {\it cross-eyed diagram twist} 
with respect to the 5-edge $\{x,y\}$. 

Let $(W,S)$ be a chordal Coxeter system of finite rank 
with just one bad edge $\{x,y\}$, and 
let ${\cal B}$ be the set of bad irreducible simplices of $(W,S)$. 
For each $A\in {\cal B}$, 
let $\beta_A:\langle A\rangle \to \langle A\rangle$ 
be an automorphism that is not inner and that maps 
each element of $A$ to a conjugate of itself in $\langle A\rangle$. 
If $S$ has no bad separators,  
let $\alpha$ be an automorphism of $W$ as in 
the automorphic angle deformation theorem, and let $S' = \alpha(S)$. 
If $S$ has a bad separator, 
let $S'$ be a set of Coxeter generators of $W$ 
as in the  5-edge angle deformation theorem. 
Then for each $A \in {\cal B}$, 
there is an element $w_A$ of $W$ such that $w_A\beta_A(A)w_A^{-1} \subset S'$  
and for each good irreducible simplex $A$ of $(W,S)$, 
there is an element $w_A$ of $W$ such that $w_AAw_A^{-1} \subset S'$. 
We say that the Coxeter system $(W,S')$ is obtained 
from the Coxeter system $(W,S)$ by an {\it elementary angle deformation}. 

\begin{theorem} 
{\rm (The Angle Deformation Theorem)}
Let $(W,S)$ be a chordal Coxeter system of finite rank, 
and let ${\cal B}$ be a set of bad irreducible simplices for $(W,S)$. 
For each $A\in {\cal B}$, 
let $\beta_A:\langle A\rangle \to \langle A\rangle$ 
be an automorphism that is not inner and that maps 
each element of $A$ to a conjugate of itself in $\langle A\rangle$. 
Then $W$ has a set of Coxeter generators $S'$ 
such that for each $A \in {\cal B}$, 
there is an element $w_A$ of $W$ such that $w_A\beta_A(A)w_A^{-1} \subset S'$  
and for each good irreducible simplex $A$ of $(W,S)$, 
there is an element $w_A$ of $W$ such that $w_AAw_A^{-1} \subset S'$. 
Moreover $(W,S')$ can be obtained from $(W,S)$ by a finite sequence 
of elementary angle deformations. 
\end{theorem}
\begin{proof}
The theorem follows easily by induction on the number of bad edges  
from the automorphic angle deformation theorem 
and the 5-edge angle deformation theorem. 
\end{proof}

\section{Conjugating Visual Subgroups}  

In this section, we prove some lemmas about conjugating visual subgroups 
of a Coxeter system. 

\begin{theorem} 
{\rm (Richardson \cite{Richardson}, Theorem A)} 
Let $(W,S)$ be a Coxeter system and let $w$ be an element of order two in $W$. 
Then there is a subset $I$ of $S$ such that $\langle I\rangle$ is finite, 
$w$ is conjugate to the maximal length element $w_I$ of $\langle I\rangle$, 
and $w_I$ is central in $\langle I \rangle$; moreover $I$ is unique 
up to conjugation in $W$. 
\end{theorem}

Let $(W,S)$ be a Coxeter system. 
An element $r$ of $W$ is called a {\it reflection} of $(W,S)$ 
if $r$ is conjugate in $W$ to an element of $S$. 
The set of all reflection of $(W,S)$ is denoted by $S^W$. 

\begin{lemma} 
Let $(W,S)$ be a Coxeter system, and 
let $r$ be a reflection of $(W,S)$. 
If $A \subset S$ and $r\in \langle A\rangle$, 
then $r$ is a reflection of $(\langle A\rangle, A)$. 
\end{lemma}
\begin{proof}
By Richardson's theorem, there is a subset $I$ of $A$ such that $\langle I\rangle$ is finite, 
$r$ is conjugate in $\langle A\rangle$ to the maximal length element $w_I$ of $\langle A\rangle$ 
and $w_I$ is central in $\langle I\rangle$. 
Now $r$ is conjugate in $W$ to an element $s$ of $S$. 
Hence $I$ is conjugate to $\{s\}$ in $W$ by Richardson's theorem. 
Therefore $I = \{a\}$ for some $a\in A$. 
Thus $r$ is a reflection of $(\langle A\rangle, A)$. 
\end{proof}

\begin{lemma} 
{\rm (Mihalik, Ratcliffe, and Tschantz \cite{M-R-T}, Lemma 6.11)}
Let $(W,S)$ be a Coxeter system with $A,B \subset S$ such that 
$\langle A\rangle$ is finite and irreducible. 
If $\langle A\rangle$ is conjugate to $\langle B\rangle$ in $W$ 
and $A$ is neither of type ${\bf A}_n$, for some $n$,  
nor of type ${\bf B}_5$, then $A=B$. 
\end{lemma}

Let $(W,S)$ be a Coxeter system, and 
let $w$ be an element of a Coxeter group $W$. 
We denote the inner automorphism $u \mapsto wuw^{-1}$ of $W$ by $w_\ast$. 
A {\it graph automorphism} of $(W,S)$ is an automorphism 
$\gamma$ of $(W,S)$. 
An {\it inner by graph automorphism} of $(W,S)$ 
is an automorphism $\alpha$ of $W$ of the form $\alpha = \iota\gamma$ 
where $\iota$ is an inner automorphism of $W$ and $\gamma$ 
is a graph automorphism of $(W,S)$. 
The set of all inner by graph automorphisms of $(W,S)$ 
forms a subgroup of ${\rm Aut}(W)$ that contains ${\rm Inn}(W)$ 
as a normal subgroup. 

Suppose $A\subset S$.  If $\langle A\rangle$ is finite, we denote 
the longest element of $\langle A\rangle$ by $\ell_A$. 
Suppose $s\in S-A$. 
Let $K\subset S$ be the irreducible component of $A\cup\{s\}$ containing $s$. 
We say that $s$ is $A$-{\it admissible} if $\langle K\rangle$ is finite. 
If $s$ is $A$-admissible, define $\nu(s,A) =\ell_K\ell_{K-\{s\}}$. 

\begin{lemma} 
Let $(W,S)$ be a Coxeter system, 
and let $A \subset S$ be an irreducible simplex. 
Let $w \in W$ be such that $w\langle A\rangle w^{-1} = \langle A\rangle$. 
Then $w_\ast$ restricts to an inner by graph automorphism on $\langle A\rangle$. 
Moreover if $A$ is of type ${\bf G}_3$, ${\bf G}_4$, or ${\bf D}_2(k)$ with $k\geq 5$, 
then $w_\ast$ restricts to an inner automorphism on $\langle A\rangle$. 
\end{lemma}
\begin{proof}
By Prop. 2.1 of Brink and Howlett \cite{B-H}, we have that $w=uv$ 
with $u\in\langle A\rangle$ and $v\in W$ with $vE_A = E_A$. 
See Prop. 4.8 of \cite{M-R-T}, for an explanation of the condition $vE_A = E_A$. 
If $\langle A\rangle$ is infinite, then $s \in S-A$ is $A$-admissible if and only if 
$s$ centralizes $A$.  Therefore $v$ centralizes $A$ by Prop. 4.10 of \cite{M-R-T}. 
Hence we may assume that $A$ is spherical. 
Now $w_\ast: \langle A\rangle \to \langle A\rangle$ preserves reflections by Lemma 7.2. 
According to Franzsen's thesis \cite{Franzsen}, 
the only irreducible spherical types that have reflection preserving automorphisms 
that are not inner by graph are ${\bf G}_3$, ${\bf G}_4$, or ${\bf D}_2(k)$ with $k\geq 5$. 
Hence we may assume that $A$ is of type ${\bf G}_3$, ${\bf G}_4$, or ${\bf D}_2(k)$ with $k\geq 5$. 
By Lemma 7.3 and Prop. 4.10 of \cite{M-R-T}, 
we have that $v$ is a product of elements of the form $\nu(s,A)$ 
such that $s$ is $A$-admissible. Suppose $s \in S-A$ is $A$-admissible. 
If $s$ centralizes $A$, then $\nu(s,A) = s$ centralizes $A$. 
If $s$ does not centralizes $A$, then $A$ is of type ${\bf D}_2(5)$ or ${\bf G}_3$ 
and $A\cup\{s\}$ is of type ${\bf G}_3$ or ${\bf G}_4$, respectively, and so
$\nu(s,A) = \ell_{A\cup\{s\}}\ell_A$ and $\ell_{A\cup\{s\}}$ centralizes $A$.  
Hence if $A$ is of type ${\bf G}_3$, ${\bf G}_4$, or ${\bf D}_2(k)$ with $k\geq 5$, 
then $w_\ast$ restricts to an inner automorphism on $\langle A\rangle$. 
\end{proof}

\begin{lemma} 
{\rm (Brady, McCammond, M\"uhlherr, Neumann \cite{B-M-M-N}, Lemma 3.7)}
If $S$ and $S'$ are two Coxeter generating sets for a Coxeter group $W$ 
and $S' \subset S^W$, then $(S')^W = S^W$. 
\end{lemma}

\begin{lemma} 
Let $(W,S)$ be a Coxeter system of finite rank, 
and let $S'$ be another set of Coxeter generators of $W$ such that $S' \subset S^W$. 
Let $A \subset S$ be an irreducible simplex, 
and let $B \subset S$ be a maximal irreducible simplex containing $A$.  
Then there is an irreducible simplex $A' \subset S'$ 
and an element $u$ of $W$ such that $\langle A\rangle = u\langle A'\rangle u^{-1}$, 
and there is an isomorphism $\alpha_A: (\langle A\rangle,A) \to (\langle A'\rangle, A')$ 
such that $u_\ast\alpha_A: \langle A\rangle \to \langle A\rangle$ maps each element of $A$ 
to a conjugate of itself in $\langle A\rangle$,  
and either $u_\ast\alpha_A$ is an inner by graph automorphism of $(\langle A\rangle,A)$, 
or $A$ and $B$ are of type ${\bf G}_3$, ${\bf G}_4$, or ${\bf D}_2(k)$ with $k\geq 5$.    
\end{lemma}
\begin{proof}
This is clear if $|A| = 1$, since $S \subset (S')^W$ by Lemma 7.5. 
Hence we may assume that $|A| > 1$. 
Let $C \subset S$ be a maximal simplex containing $B$. 
By Prop. 5.10 of \cite{M-R-T}, there is a maximal simplex $C'\subset S'$ 
and an element $w$ of $W$ such that $\langle C'\rangle = w\langle C\rangle w^{-1}$. 
The isomorphism $w_\ast:\langle C\rangle \to \langle C'\rangle$ 
preserves reflections by Lemma 7.2. 
By Lemma 14 of Franzsen and Howlett \cite{F-H}, 
$w_\ast$ maps $B$ onto a component $B'$ of $C'$, 
and so $w\langle B\rangle w^{-1} = \langle B'\rangle$. 

If $\langle B\rangle$ is finite, then $B$ and $B'$ have the same type, 
since $w_\ast:\langle B\rangle \to \langle B'\rangle$ is an isomorphism 
and $B$ and $B'$ are both irreducible. 
If $\langle B\rangle$ is infinite, then $B$ is strongly rigid by Theorem 2.2. 
Hence $B$ and $B'$ have the same type. 
Thus there is an isomorphism $\alpha: (\langle B\rangle,B) \to (\langle B'\rangle,B)$. 
The isomorphism $w_\ast^{-1}\alpha:\langle B\rangle \to \langle B\rangle$ 
preserves reflections by Lemma 7.2. 

Assume first that $w_\ast^{-1}\alpha = v_\ast\gamma$ 
where $v\in \langle B\rangle$ and $\gamma$ is a graph automorphism of $(\langle B\rangle,B)$. 
Then by replacing $\alpha$ by $\alpha\gamma$, we may assume that $w_\ast^{-1}\alpha = v_\ast$. 
Let $A' = \alpha(A)$ and $u = v^{-1}w^{-1}$. 
Then $\langle A'\rangle  = u^{-1}\langle A\rangle u$, 
and so $\langle A\rangle = u\langle A'\rangle u^{-1}$. 
Let $\alpha_A:(\langle A\rangle, A) \to (\langle A'\rangle, A')$ 
be the restriction of $\alpha$. 
Then $u_\ast\alpha_A: \langle A\rangle \to \langle A\rangle$ is the identity map, 
since $u_\ast\alpha: \langle B\rangle \to \langle B\rangle$ is the identity map.

Now assume $w_\ast^{-1}\alpha:\langle B\rangle \to \langle B\rangle$ is not inner by graph. 
If $\langle B\rangle$ is infinite, then every automorphism of $\langle B\rangle$ 
is inner by graph, since $\langle B\rangle$ is strongly rigid. 
Hence $B$ is spherical and 
$B$ is of type ${\bf G}_3$, ${\bf G}_4$, or ${\bf D}_2(k)$ with $k\geq 5$. 
By Prop. 32 of Franzsen and Howlett \cite{F-H}, 
there is a $v\in \langle B\rangle$ such that 
$w_\ast^{-1}\alpha\langle A\rangle = v\langle A\rangle v^{-1}$. 
Let $A' = \alpha(A)$ and $u = v^{-1}w^{-1}$. 
Then $\langle A'\rangle  = u^{-1}\langle A\rangle u$, 
and so $\langle A\rangle = u\langle A'\rangle u^{-1}$. 
Let $\alpha_A:(\langle A\rangle, A) \to (\langle A'\rangle, A')$ 
be the restriction of $\alpha$. 
If $u_\ast\alpha_A:\langle A\rangle \to \langle A\rangle$ does not map each 
element of $A$ to a conjugate of itself in $\langle A\rangle$, then 
$A = B$ and $B$ is of type ${\bf D}_2(k)$ with $k$ even, 
in which case, we replace $\alpha$ by $\alpha\gamma$ 
where $\gamma$ is the nontrivial graph automorphism of $(\langle B\rangle,B)$. 
Then  $u_\ast\alpha_A:\langle A\rangle \to \langle A\rangle$ maps each 
element of $A$ to a conjugate of itself in $\langle A\rangle$. 
Finally, if $u_\ast\alpha_A:\langle A\rangle \to \langle A\rangle$ is not inner by graph, 
then $A$ is of type ${\bf G}_3$, ${\bf G}_4$, or ${\bf D}_2(k)$ with $k\geq 5$. 
\end{proof} 

\begin{lemma} 
Let $(W,S)$ be a Coxeter system of finite rank, 
and let $S'$ be another set of Coxeter generators of $W$ such that $S' \subset S^W$. 
Let $A \subset S$ be an irreducible simplex, 
and let $A' \subset S'$ and  $\ov A\hbox{}' \subset S'$ be irreducible simplices, 
and let $u,\ov u\in W$ such that 
$\langle A\rangle = u\langle A'\rangle u^{-1} = \ov u\langle \ov A\hbox{}'\rangle \ov u\hbox{}^{-1}$. 
Let $\alpha:(\langle A\rangle,A) \to (\langle A'\rangle,A')$ be an isomorphism, 
and let $\ov\alpha: (\langle A\rangle,A) \to (\langle \ov A\hbox{}'\rangle,A')$ 
be an isomorphism. 
Then $u_\ast\alpha: \langle A\rangle \to \langle A\rangle$ is an inner by graph 
automorphism of $(\langle A\rangle, A)$ 
if and only if $\ov u_\ast\ov\alpha: \langle A\rangle \to \langle A\rangle$ is an inner by graph 
automorphism of $(\langle A\rangle, A)$.  
\end{lemma}
\begin{proof}
Now $u^{-1}\ov u\langle \ov A\hbox{}'\rangle\ov u\hbox{}^{-1}u = \langle A'\rangle$. 
Let $v$ be the shortest element of the double coset 
$\langle A'\rangle u^{-1}\ov u\langle \ov A\hbox{}'\rangle$. 
Then $v\ov A\hbox{}'v^{-1} = A'$ by Lemma 4.3 of \cite{M-R-T}. 
Now $\ov u_\ast\ov\alpha$ is inner by graph if and only if $(\ov uv^{-1})_\ast(v_\ast\ov\alpha)$ 
is inner by graph. 
Hence by replacing $\ov A\hbox{}'$ by $A'$, $\ov\alpha$ by $v_\ast\ov\alpha$, and $\ov u$ 
by $\ov uv^{-1}$, we may assume that $\ov A\hbox{}'=A'$. 
Now $\ov u_\ast\ov\alpha = (\ov uu^{-1})_\ast(u_\ast\alpha)(\alpha^{-1}\ov\alpha)$. 
Observe that $\alpha^{-1}\ov\alpha$ is a graph automorphism of $\langle A\rangle$   
and $(\ov uu^{-1})_\ast$ restricts to an inner by graph automorphism of $\langle A\rangle$ 
by Lemma 7.4.  
Therefore $u_\ast\alpha: \langle A\rangle \to \langle A\rangle$ is inner by graph 
if and only if $\ov u_\ast\ov\alpha: \langle A\rangle \to \langle A\rangle$ is inner by graph.
\end{proof}

\section{The Sharp Angle Theorem} 

Let $(W,S)$ be a Coxeter system, and let $S'$ be another set of Coxeter generators of $W$ 
such that $S' \subset S^W$. 
Then $S$ is said to be {\it sharp-angled} with respect to $S'$ if for each pair 
$s,t\in S$ such that $2<m(s,t)<\infty$, 
there is a $w\in W$ such that $w\{s,t\}w^{-1} \subset S'$. 

\begin{lemma}  
Let $(W,S)$ be a Coxeter system, and let $S'$ be another set of Coxeter generators of $W$ 
such that $S' \subset S^W$. 
Then the following are equivalent:
\begin{enumerate}
\item The set $S$ is sharp-angled with respect to $S'$.
\item If $A\subset S$ is an irreducible simplex of rank 2, and  
$A'\subset S'$ is an irreducible simplex of rank 2, and 
$u \in W$ such that $\langle A\rangle = u\langle A'\rangle u^{-1}$, 
and $\alpha:(\langle A\rangle,A) \to (\langle A'\rangle, A')$ is an isomorphism, 
then $u_\ast\alpha:\langle A\rangle \to \langle A\rangle$ is an inner by graph automorphism
of $(\langle A\rangle, A)$.
\item If $A\subset S$ is an irreducible simplex, and 
$A'\subset S'$ is an irreducible simplex, and 
$u \in W$ such that $\langle A\rangle = u\langle A'\rangle u^{-1}$, 
and $\alpha:(\langle A\rangle,A) \to (\langle A'\rangle, A')$ is an isomorphism, 
then $u_\ast\alpha:\langle A\rangle \to \langle A\rangle$ is an inner by graph automorphism 
of $(\langle A\rangle, A)$.
\item If $A\subset S$ is an irreducible simplex, 
then there is $w\in W$ such that $wAw^{-1} \subset S'$.
\end{enumerate}
\end{lemma}
\begin{proof}
Suppose $S$ is sharp-angled with respect to $S'$. 
Let $A\subset S$ be an irreducible simplex of rank 2, let  
$A'\subset S'$ be an irreducible simplex of rank 2, let 
$u \in W$ be such that $\langle A\rangle = u\langle A'\rangle u^{-1}$, 
and let $\alpha:(\langle A\rangle,A) \to (\langle A'\rangle, A')$ be an isomorphism. 
Then there exists $w\in W$ such that $wAw^{-1} \subset S'$. 
Let $\ov A\hbox{}'=wAw^{-1}$. 
Define $\ov\alpha:\langle A\rangle \to\langle\ov A\hbox{}'\rangle$ by $\ov\alpha(a)=waw^{-1}$. 
Then $w_\ast^{-1}\ov\alpha:\langle A\rangle \to\langle A\rangle$ is the identity map. 
Hence $u_\ast\alpha:\langle A\rangle \to \langle A\rangle$ is inner by graph by Lemma 7.7. 
Hence (1) implies (2). 

Assume (2) holds and $A\subset S$ is an irreducible simplex, and 
$A'\subset S'$ is an irreducible simplex, and 
$u \in W$ such that $\langle A\rangle = u\langle A'\rangle u^{-1}$, 
and $\alpha:(\langle A\rangle,A) \to (\langle A'\rangle, A')$ is an isomorphism. 
As $S \subset (S')^W$ and (2) holds, we may assume that ${\rm rank}(A) > 2$. 
On the contrary, assume $u_\ast\alpha:\langle A\rangle \to \langle A\rangle$ 
is not inner by graph. 
As infinite complete Coxeter systems are strongly rigid, $\langle A\rangle$ must be finite. 
The automorphism $u_\ast\alpha$ preserves reflections, since $S'\subset S^W$. 
Therefore $A$ is of type ${\bf G}_3$ or ${\bf G}_4$. 
Let $B$ be the subset of $A$ of type ${\bf D}_2(5)$ and let $B'$ be the subset of $A'$ 
of type ${\bf D}_2(5)$. 
Let $\ov\alpha:(\langle B\rangle,B)\to(\langle B'\rangle,B')$ be the restriction 
of $\alpha:(\langle A\rangle,A)\to (\langle A'\rangle,A')$. 
By Prop. 32 of Franzsen and Howlett \cite{F-H}, 
there is $v\in \langle A\rangle$ such that $vu\langle B'\rangle u^{-1}v^{-1} = \langle B\rangle$. 
By a computer calculation of all possibilities, 
we find that $(vu)_\ast\ov\alpha:\langle B\rangle \to \langle B\rangle$ 
is not inner by graph, which contradicts (2).  Therefore (2) implies (3).

Now assume (3) holds.  Let $A\subset S$ be an irreducible simplex. 
By Lemma 7.6, there is an irreducible simplex $A'\subset S'$ and 
an element $u$ of $W$ such that $\langle A\rangle = u\langle A'\rangle u^{-1}$. 
Let $\alpha:(\langle A\rangle,A)\to (\langle A'\rangle,A')$ be an isomorphism. 
Then $u_\ast\alpha:\langle A\rangle \to \langle A\rangle$ is inner by graph. 
Hence there is a $v\in\langle A\rangle$ and a graph automorphism $\gamma$ of $(\langle A\rangle,A)$ 
such that $u_\ast\alpha= v_\ast\gamma$.  Hence $u\alpha(A)u^{-1} = v\gamma(A)v^{-1}$, 
and so $uA'u^{-1} = vAv^{-1}$. 
Therefore $u^{-1}vAv^{-1}u = A'\subset S'$. 
Thus (3) implies (4). 

Now (4) for ${\rm rank}(A) = 2$ implies that $S$ is sharp-angled with respect to $S'$. 
Thus (4) implies (1). 
\end{proof}

\begin{theorem} 
{\rm (The Sharp Angle Theorem)} 
Let $(W,S)$ be a chordal Coxeter system of finite rank, 
and let $S'$ be another set of Coxeter generators 
of $W$ such that $S'\subset S^W$. 
Then there is a set of Coxeter generators $S_\ast$ of $W$ 
such that $S_\ast \subset S^W$, the set $S_\ast$ is sharp-angled with respect to $S'$, 
and $S_\ast$ can be obtained from $S$ by a finite sequence of elementary angle deformations. 
\end{theorem}
\begin{proof}
Let $A\subset S$ be an irreducible simplex. 
By Lemma 7.6, there is an irreducible simplex $A'\subset S'$ 
and an element $u_A$ of $W$ such that $\langle A\rangle = u_A\langle A'\rangle u_A^{-1}$, 
and there is an isomorphism $\alpha_A: (\langle A\rangle,A) \to (\langle A'\rangle, A')$ 
such that $(u_A)_\ast\alpha_A: \langle A\rangle \to \langle A\rangle$ maps each element of $A$ 
to a conjugate of itself in $\langle A\rangle$,  
and either $(u_A)_\ast\alpha_A$ is an inner by graph automorphism of $(\langle A\rangle,A)$, 
or $A$ is of type ${\bf G}_3$, ${\bf G}_4$, or ${\bf D}_2(k)$ with $k\geq 5$. 
By Lemma 7.7, the condition that $(u_A)_\ast\alpha_A$ is inner by graph depends only on $A$. 
Let ${\cal B}$ be the set of all irreducible simplices $A$ of $(W,S)$ 
such that $(u_A)_\ast\alpha_A$ is not inner by graph. 
Then ${\cal B}$ is a set of bad irreducible simplices for $(W,S)$ by Lemma 7.6. 

For each $A \in {\cal B}$, let $\beta_A = (u_A)_\ast\alpha_A$. 
By the angle deformation theorem, there is a set of Coxeter generators $S_\ast$ for $W$ 
such that for each $A \in {\cal B}$, 
there is a $w_A\in W$ such that $w_A\beta_A(A)w_A^{-1} \subset S_\ast$  
and for each good irreducible simplex $A$ of $(W,S)$, 
there is a $w_A\in W$ such that $w_AAw_A^{-1} \subset S_\ast$. 
Moreover $(W,S_\ast)$ can be obtained from $(W,S)$ by a finite sequence 
of elementary angle deformations. 

As singleton simplices of $(W,S)$ are good, we have that $S \subset S_\ast^W$. 
By Lemma 7.5, we have that $S_\ast^W = S^W = (S')^W$. 
Let $s_\ast, t_\ast$ be elements of $S_\ast$ such that $2 < m(s_\ast,t_\ast) < \infty$, 
and let $A_\ast =\{s_\ast,t_\ast\}$.  
We need to show that there is a $w\in W$ such that $wA_\ast w^{-1} \subset S'$. 
As every automorphism of ${\bf A}_2$ and ${\bf C}_2$ is inner by graph, 
we may assume that $m(s_\ast,t_\ast) \geq 5$ by Lemma 8.1.  
Then there is a unique irreducible simplex $A = \{s,t\}$ of $(W,S)$ 
such that $\langle A\rangle$ is conjugate to $\langle A_\ast\rangle$ in $W$ by Lemmas 7.3 and 7.6. 
Let $A'=\alpha_A(A)$. 

First assume that $(u_A)_\ast\alpha_A = v_\ast\gamma$ 
where $v\in \langle A\rangle$ and $\gamma$ is a graph automorphism of $(\langle A\rangle,A)$. 
Then $u_A A'u_A^{-1} = vAv^{-1}$. 
Now $w_A A w_A^{-1} \subset S_\ast$ implies that $w_A A w_A^{-1} = A_\ast$ by Lemma 7.3. 
Let $w=u_A^{-1}vw_A^{-1}$. Then 
$$wA_\ast w^{-1} = u_A^{-1}vw_A^{-1}A_\ast w_Av^{-1}u_A^{-1} = u_A^{-1}vAv^{-1}u_A = A' \subset S'.$$

Now assume that $(u_A)_\ast\alpha_A$ is not inner by graph. 
Then we have that 
$w_A\beta_A(A)w_A^{-1} \subset S_\ast,$
and so $w_A\beta_A(A)w_A^{-1} = A_\ast$ by Lemma 7.3. 
Hence $w_Au_AA' u_A^{-1}w_A^{-1} = A_\ast$. 
Let $w=u_A^{-1}w_A^{-1}$. 
Then $wA_\ast w^{-1} = A' \subset S'$. 
Thus $S_\ast$ is sharp-angled with respect to $S'$. 
\end{proof}

\section{The Twist Equivalence Theorem} 

\begin{lemma} 
Let $(W,S)$ be a Coxeter system of finite rank and let $S'$ 
be another set of Coxeter generators of $W$ such that $S' \subset S^W$ 
and $S$ is sharp-angled with respect to $S'$. 
Let $B\subset S$ and $B'\subset S'$ such that $\langle B\rangle = \langle B'\rangle$. 
Then $B'\subset B^{\langle B\rangle}$ and $B$ is sharp-angled with respect to $B'$. 
\end{lemma}
\begin{proof}
By Lemma 7.2, we have that $B'\subset B^{\langle B\rangle}$. 
Let $A\subset B$ be an irreducible simplex of rank 2. 
By Lemma 7.6, there is an irreducible simplex $A'\subset B'$ of rank 2 
and a $u\in \langle B\rangle$ such that $u\langle A'\rangle u^{-1} =\langle A\rangle$. 
Let $\beta:(\langle A\rangle,A)\to(\langle A'\rangle,A')$ be an isomorphism. 
Then $u_\ast\beta:\langle A\rangle \to\langle A\rangle$ is inner by graph by Lemma 8.1, 
since $S$ is sharp-angled with respect to $S'$. 
Therefore $B$ is sharp-angled with respect $B'$ by Lemmas 7.7 and 8.1. 
\end{proof} 

\begin{lemma} 
Let $(W,S)$ be a Coxeter system 
and let $\tau$ be the elementary twist of $(W,S)$, with respect to a 
separation $(S_1,S_0,S_2)$ of $S$, by the element $\ell$ of $\langle S_0\rangle$. 
Here $\tau: S\mapsto S_1\cup\ell S_2\ell^{-1}$. 
Let $\tau'$ be the elementary twist of $(W,S)$, with respect to the separation 
$(S_2,S_0,S_1)$ of $S$, by $\ell$. 
Here $\tau': S\mapsto \ell S_1\ell^{-1}\cup S_2$. 
Then $\ell^2=1$ and $\tau = \ell_\ast\tau'$.
\end{lemma} 

\begin{lemma} 
Let $(W,S)$ be a Coxeter system 
and let $\tau$ be the elementary twist of $(W,S)$, with respect to a 
separation $(S_1,S_0,S_2)$ of $S$, by the element $\ell$ of $\langle S_0\rangle$. 
Let $\alpha$ be an automorphism of $W$ of order 2. 
Let $\alpha_\ast(\tau)$ be the elementary twist of $(W,\alpha(S))$, 
with respect to the separation $(\alpha(S_1),\alpha(S_0),\alpha(S_2))$ of $\alpha(S)$, 
by the element $\alpha(\ell)$ of $\langle\alpha(S_0)\rangle$. 
Here 
$$\alpha_\ast(\tau):\alpha(S)\mapsto \alpha(S_1)\cup\alpha(\ell)\alpha(S_2)\alpha(\ell)^{-1}.$$
Then $\tau\alpha = \alpha\alpha_\ast(\tau): \alpha(S)\mapsto S_1\cup\ell S_2\ell^{-1}$. 
\end{lemma}

\begin{theorem} 
{\rm (The Twist Equivalence Theorem)}
Let $(W,S)$ be a chordal Coxeter system of finite rank, 
and let $S'$ be another set of Coxeter generators for $W$ 
such that $S'\subset S^W$ and $S$ is sharp-angled with respect to $S'$. 
Then $(W,S)$ is twist equivalent to $(W,S')$. 
\end{theorem}
\begin{proof}
The proof is by induction on $|S|$. 
The theorem is clear if $|S| =1$, 
so assume $|S| > 1$ and the theorem 
is true for all chordal Coxeter systems 
of rank less than $|S|$. 
Assume first that $(W,S)$ is complete. 
Then $(W,S')$ is complete by Prop. 5.10 of \cite{M-R-T}. 
Let 
$$(W,S) = (W_1,S_1)\times\cdots\times(W_n,S_n)$$
be the factorization of $(W,S)$ into irreducible factors, and let 
$$(W,S') = (W_1',S_1')\times\cdots\times(W_m',S_m')$$
be the factorization of $(W,S')$ into irreducible factors. 
By Lemma 14 of Franzsen and Howlett \cite{F-H}, 
$m=n$ and by reindexing, we may assume that $W_i'=W_i$ for each $i=1,\ldots,n$. 
As $S$ is sharp-angled with respect to $S'$, there is a $w_i\in W$ such that 
$w_iS_iw_i^{-1}\subset S'$ for each $i$ by Lemma 8.1. 
As the $j$th component of $w_i$, for $j\neq i$, centralizes $W_i$, 
we may assume that $w_i\in W_i$. 
Then $w_iS_iw_i^{-1} = S_i'$ for each $i$. 
Let $w=w_1\cdots w_n$. 
Then $wSw^{-1} = S'$. 
Therefore $(W,S)$ is twist equivalent to $(W,S')$. 

Now assume that $(W,S)$ is incomplete. 
Then $\Gamma(W,S)$ has a cut set. 
Let $C \subset S$ be a minimal cut set of $\Gamma(W,S)$. 
By Theorem 1 of Dirac \cite{Dirac}, 
we have that $C$ is a simplex.  
Let $(A,C,B)$ be the separation of $S$ determined by $C$. 
Then $W = \langle A\rangle \ast_{\langle C\rangle} \langle B\rangle$ 
is a nontrivial splitting. 
By the Decomposition Matching Theorem, 
$(W,S)$ and $(W,S')$ are twist equivalent 
to Coxeter systems $(W,R)$ and $(W,R')$, respectively, 
such that there exists a nontrivial
reduced visual graph of groups decomposition $\Psi$ of $(W,R)$
and a nontrivial reduced visual graph of groups decomposition
$\Psi'$ of $(W,R')$ having the same graphs and the same vertex
and edge groups and all edge groups equal and a subgroup of a
conjugate of $\langle C \rangle$.
The Coxeter systems $(W,R)$ and $(W,S)$ are twist equivalent, 
and so $|R| = |S|$ and $(W,R)$ is chordal by Lemma 3.2. 
Now $R'\subset R^W$ and $R$ is sharp-angled with respect to $R'$, 
since $R'$ is twist equivalent to $S'$. 

Let $\{(W_i,R_i)\}_{i=1}^k$ be the Coxeter systems of the vertex groups of $\Psi$, 
and let $(W_0,R_0)$ be the Coxeter system of the edge group of $\Psi$. 
Then $k\geq 2$, and $R =\cup_{i=1}^k R_i$, and $\cap_{i=1}^k R_i = R_0$, 
and $R_i - R_0 \neq \emptyset$ for each $i > 0$, 
and $m(a,b) = \infty$ for each $a \in R_i-R_0$ and $b \in R_j-R_0$ with $i\neq j$. 
From the proof of Lemma 8.4 of \cite{M-R-T}, we deduce that $R_0$ 
is conjugate to a subset of $S$. 
By Lemma 4.3 of \cite{M-R-T}, we have that $R_0$ is conjugate to a subset of $C$, 
and so $R_0$ is a simplex. 
Let $\{(W_i',R_i')\}_{i=1}^k$ be the Coxeter systems of the vertex groups of $\Psi'$ 
indexed so that $W_i' = W_i$ for each $i$, 
and let $(W_0',R_0')$ be the Coxeter system of the edge group of $\Psi'$. 
Then $W_0' = W_0$, and $R' =\cup_{i=1}^k R_i'$, and $\cap_{i=1}^k R_i' = R_0'$, 
and $R_i' - R_0' \neq \emptyset$ for each $i > 0$, 
and $m(a',b') = \infty$ for each $a' \in R_i'-R_0'$ and $b' \in R_j'-R_0'$ with $i\neq j$. 

The Coxeter system $(W_i,R_i)$ is chordal and $|R_i| < |R|$ for each $i$. 
By Lemma 9.1, we have  $R_i'\subset R_i^{W_i}$ and  
$R_i$ is sharp-angled with respect to $R_i'$ for each $i$. 
Hence by the induction hypothesis, $(W_i,R_i)$ is twist equivalent to $(W_i,R_i')$ for each $i$. 
As $R_0$ is a simplex, there is an element $w_0$ of $W_0$ such that $w_0R_0w_0^{-1} = R_0'$. 
By conjugating $W$ by $w_0$, we may assume that $R_0=R_0'$. 

Now there is a sequence $\tau_1,\ldots,\tau_m$ of elementary twists 
that transforms $R_1$ to $R_1'$. 
Suppose $\tau_1$ is the elementary twist of $(W_1,R_1)$, 
with respect to a separation $(S_1,S_0,S_2)$ of $R_1$, 
by an element $\ell_1$ of $\langle S_0\rangle$. 
Here $\tau_1:R_1\mapsto S_1\cup\ell_1S_2\ell_1^{-1}$. 
As $R_0$ is a simplex, either $R_0\subset S_1$ or $R_0\subset S_2$. 
If $R_0\subset S_1$, let $\tau_1'=\tau_1$. 
If $R_0\subset S_2$, replace $\tau_1$ by $(\ell_1)_\ast\tau_1'$ 
where $\tau_1'$ is the elementary twist of $(W_1,R_1)$ by $\ell_1$, with respect to 
the separation $(S_2,S_0,S_1)$ of $R_1$. 
Here $\tau_1': R_1\mapsto S_2\cup\ell_1S_1\ell_1^{-1}$. 
By Lemma 9.3, we can move $(\ell_1)_\ast$ past $\tau_2,\ldots, \tau_m$ to get 
a sequence of elementary twists $\tau_1',\ov\tau_2,\ldots,\ov\tau_m$ 
followed by $(\ell_1)_\ast$ that transforms $R_1$ to $R_1'$. 
Repeating this procedure, we obtain a sequence of elementary twists 
$\tau_1',\tau_2',\ldots,\tau_m'$ followed by an inner automorphism $(w_1)_\ast$  
with $w_1\in W_1$ that transforms $R_1$ to $R_1'$ 
such that $R_0$ is part of the left (invariant) side of the separation of each twist. 
Hence the elementary twists $\tau_1',\tau_2',\ldots,\tau_m'$ extend to all of $W$, 
with $R_i$ for $i>1$, part of the left (invariant) side of the separation of each twist. 
Let the sequence $\tau_1',\tau_2',\ldots,\tau_m'$  transform $R_1$ to $R_1^\ast$. 
Then $R_0\subset R_1^\ast$. 
As $w_1R_1^\ast w_1^{-1} = R_1'$ and $R_0\subset R_1'$. 
there is a subset $R_0^\ast$ of $R_1^\ast$ such that $w_1R_0^\ast w_1^{-1} = R_0$. 
Hence the generalized twist of $(W,R_1^\ast\cup\cup_{i>1}R_i)$ by $w_1\in W_1$, 
with respect to the separation $(R_1^\ast,R_0,\cup_{i>1}R_i)$ of $R_1^\ast\cup\cup_{i>1}R_i$, 
transforms $R_1^\ast\cup\cup_{i>1}R_i$ to $R_1'\cup\cup_{i>1}R_i$. 
Likewise $(W,R_1'\cup\cup_{i>1}R_i)$ is twist equivalent to $(W,R_1'\cup R_2'\cup\cup_{i>2}R_i)$, 
and by induction $(W,R)$ is twist equivalent to $(W,R')$. 
Thus $(W,S)$ is twist equivalent to $(W,S')$. 
\end{proof}

\section{The Isomorphism Theorem}

Let $(W,S)$ be a Coxeter system. 
A {\it basic subgroup} of $(W,S)$ is a noncyclic, maximal, finite, irreducible, 
visual subgroup of $(W,S)$. 
A ${\it base}$ of $(W,S)$ is a subset $B$ of $S$ such that 
$\langle B\rangle$ is a basic subgroup of $(W,S)$. 

Let $S'$ be another set of Coxeter generators of $W$. 
A base $B$ of $(W,S)$ is said to {\it match} a base $B'$ of $(W,S')$ 
if $[\langle B\rangle,\langle B\rangle]$ is conjugate to 
$[\langle B'\rangle,\langle B'\rangle]$ in $W$.  
The next two theorems were proved in \S 5 of \cite{M-R-T}.

\begin{theorem} 
Let $B$ be a base of $(W,S)$ of type ${\bf C}_{2q+1}$ for some $q\geq 1$, and   
let $a,b,c$ be the elements of $B$ such that $m(a,b)=4$ and $m(b,c)=3$. 
Suppose that $m(s,t)=2$ for all $(s,t)\in (S-B)\times B$ such that $m(s,a)<\infty$. 
Let $d =aba$, and let $z$ be the longest element of $\langle B\rangle$. 
Let $S'=(S-\{a\})\cup\{d,z\}$ and $B'=(B-\{a\})\cup\{d\}$. 
Then $S'$ is a set of Coxeter generators for $W$ 
such that 
\begin{enumerate}
\item The set $B'$ is a base of $(W,S')$ of type ${\bf B}_{2q+1}$ that matches $B$,  
\item $m(z,t)=2$ for all $t\in B'$,
\item If $(s,t)\in (S-B)\times\{d,z\}$, then 
$m(s,t)<\infty$ if and only if $m(s,a)<\infty$, moreover  
if $m(s,t)<\infty$, then $m(s,t)=2$.
\end{enumerate}
\end{theorem}

\begin{theorem} 
Let $B=\{a,b\}$ be a base of $(W,S)$ of type ${\bf D}_2(4q+2)$ for some $q\geq 1$. 
Suppose that if $s\in S-B$ and $m(s,a)<\infty$, then $m(s,a) = m(s,b) = 2$. 
Let $c=aba$ and let $z$ be the longest element of $\langle B\rangle$. 
Let $S'=(S-\{a\})\cup\{c,z\}$ and $B'=\{b,c\}$. 
Then $S'$ is a set of Coxeter generators of $W$ such that 
\begin{enumerate}
\item The set $B'$ is a base of $(W,S')$ of type ${\bf D}_2(2q+1)$ that matches $B$,  
\item $m(z,b)=m(z,c) = 2$, 
\item if $(s,t)\in (S-B)\times \{c,z\}$, then $m(s,t)<\infty$ if and only if $m(s,a)<\infty$,  
moreover if $m(s,t)<\infty$, then $m(s,t) = 2$. 
\end{enumerate}
\end{theorem}

Let $(W,S)$ be a Coxeter system of finite rank. 
We say that $(W,S)$ can be {\it blown up along a base} $B$ if 
$(W,S)$ and $B$ satisfy the hypothesis of either Theorem 10.1 or 10.2. 
If $(W,S)$ can be blown up along a base $B$, 
then we can blow up $(W,S)$ to a Coxeter system $(W,S')$ as in 
the statement of Theorem 10.1 or 10.2 such that $|S'| = |S|+1$, 
the base $B$ matches a base $B'$ of $(W,S')$ 
with $|\langle B\rangle| >  |\langle B'\rangle|$,  
and each other base $C$ of $(W,S)$ is also a base of $(W,S')$. 
We say that $(W,S')$ is obtained by {\it blowing up} $(W,S)$ {\it along the base} $B$. 

By the process of blowing up along a base, 
we can effectively construct a sequence 
$S = S^{(0)}, S^{(1)},\ldots, S^{(\ell)}$ of Coxeter generators of $W$ 
such that $(W,S^{(i+1)})$ is obtained by blowing up $(W, S^{(i)})$ 
along a base for each $i = 0,\ldots, \ell-1$ and $(W,S^{(\ell)})$ cannot be blown up along a base. 
The sequence terminates since the sum of the orders 
of the basic subgroups decreases at each step of the sequence.  
By Theorem 10.1 of \cite{M-R-T}, the system $(W, S^{(\ell)})$ has maximum rank 
over all Coxeter systems for $W$.  
We call $(W, S^{(\ell)})$ an {\it expanded system} determined by $(W,S)$.  

\begin{theorem} 
{\rm (Howlett and M\"uhlherr \cite{H-M}, Theorem 1)}\ \, 
If $(W,S)$ and $(W,S')$ are expanded Coxeter systems of finite rank, 
then there is an automorphism $\alpha$ of $W$ 
such that $\alpha(S) \subset (S')^W$. 
\end{theorem}

\begin{theorem} 
{\rm (The Isomorphism Theorem)}
Let $(W,S)$ and $(W',S')$ be chordal Coxeter systems of finite rank 
and let $(W,S_\ast)$ and $(W,S'_\ast)$ be expanded 
systems determined by $(W,S)$ and $(W',S')$, respectively.
Then $W$ is isomorphic to $W'$ if and only if 
$\Gamma(W,S_\ast)$ can be deformed into a labeled diagram 
isomorphic to $\Gamma(W',S_\ast')$ by a finite sequence 
of elementary diagram twists and cross-eyed diagram twists. 
\end{theorem}
\begin{proof}
The Isomorphism Theorem follows from Theorem 10.3, 
the 5-Edge Angle Deformation Theorem, the Sharp Angle Theorem, and 
the Twist Equivalence Theorem. 
\end{proof}

Consider the chordal Coxeter system $(W,S)$ defined by the 
upper left diagram in Figure 2. 
By the Isomorphism Theorem, Figure 2 illustrates all four possible isomorphism types 
of P-diagrams for $W$. 
The diagrams in each row differ by an elementary twist along the 5-edge. 
The diagrams in each column differ by a cross-eyed twist with respect 
to the 5-edge.

\medskip

$$\mbox{
\setlength{\unitlength}{.8cm}
\begin{picture}(16,13.5)(0,0)
\thicklines
\put(4,2){\circle*{.15}}
\put(4,2){\line(0,1){3}}
\put(4,5){\circle*{.15}}
\put(2,3.5){\line(4,3){2}}
\put(2,3.5){\line(4,-3){4}}
\put(6,3.5){\line(0,-1){3}}
\put(2,3.5){\line(0,-1){3}}
\put(2,3.5){\circle*{.15}}
\put(6,3.5){\line(-4,3){2}}
\put(6,3.5){\line(-4,-3){4}}
\put(6,3.5){\circle*{.15}}
\put(2,.5){\circle*{.15}}
\put(6,.5){\circle*{.15}}
\put(1.5,1.8){2}
\put(6.3,1.8){2}
\put(5.05,2.1){2}
\put(2.75,2.1){2}
\put(2.7,4.5){3}
\put(3,0.7){3}
\put(4.7,0.7){3}
\put(4.25,3.3){5}
\put(5.2,4.5){3}
\put(12,2){\circle*{.15}}
\put(12,2){\line(0,1){3}}
\put(12,5){\circle*{.15}}
\put(10,3.5){\line(4,3){4}}
\put(14,6.5){\circle*{.15}}
\put(10,3.5){\line(4,-3){2}}
\put(10,3.5){\circle*{.15}}
\put(10,3.5){\line(0,-1){3}}
\put(10,0.5){\circle*{.15}}
\put(14,3.5){\line(-4,3){2}}
\put(14,3.5){\line(-4,-3){4}}
\put(14,3.5){\circle*{.15}}
\put(14,3.5){\line(0,1){3}}
\put(9.5,1.8){2}
\put(10.7,2.15){2}
\put(10.7,4.5){3}
\put(12.25,3.3){5}
\put(13.05,2.1){3}
\put(13.2,4.5){2}
\put(12.7,6){3}
\put(11,0.7){3}
\put(14.3,4.8){2}
\put(4,9){\circle*{.15}}
\put(4,9){\line(0,1){3}}
\put(4,12){\circle*{.15}}
\put(2,10.5){\line(4,3){2}}
\put(2,10.5){\line(4,-3){4}}
\put(6,10.5){\line(0,-1){3}}
\put(2,10.5){\line(0,-1){3}}
\put(2,10.5){\circle*{.15}}
\put(6,10.5){\line(-4,3){2}}
\put(6,10.5){\line(-4,-3){4}}
\put(6,10.5){\circle*{.15}}
\put(2,7.5){\circle*{.15}}
\put(6,7.5){\circle*{.15}}
\put(1.5,8.8){3}
\put(6.3,8.8){3}
\put(5.05,9.1){2}
\put(2.75,9.1){2}
\put(2.7,11.5){3}
\put(3,7.7){2}
\put(4.7,7.7){2}
\put(4.25,10.3){5}
\put(5.2,11.5){3}
\put(12,9){\circle*{.15}}
\put(12,9){\line(0,1){3}}
\put(12,12){\circle*{.15}}
\put(10,10.5){\line(4,3){4}}
\put(14,13.5){\circle*{.15}}
\put(10,10.5){\line(4,-3){2}}
\put(10,10.5){\circle*{.15}}
\put(10,10.5){\line(0,-1){3}}
\put(10,7.5){\circle*{.15}}
\put(14,10.5){\line(-4,3){2}}
\put(14,10.5){\line(-4,-3){4}}
\put(14,10.5){\circle*{.15}}
\put(14,10.5){\line(0,1){3}}
\put(9.5,8.8){3}
\put(10.7,9.15){2}
\put(10.7,11.5){3}
\put(12.25,10.3){5}
\put(13.05,9.1){3}
\put(13.2,11.5){2}
\put(12.7,13){2}
\put(11,7.7){2}
\put(14.3,11.8){3}

\end{picture}}$$

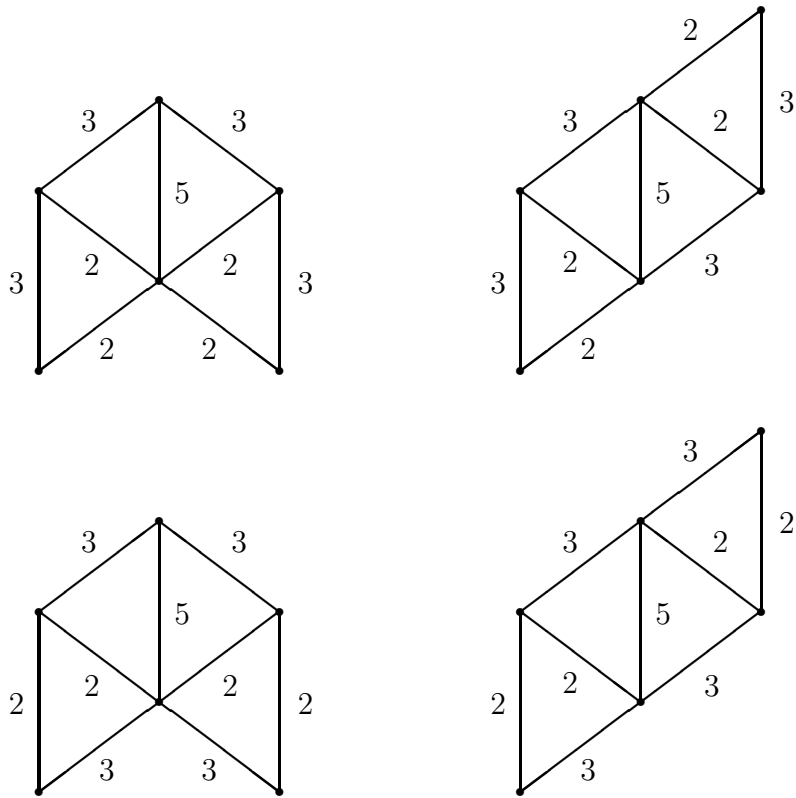
\centerline{\bf Figure 2}
\medskip

\section{A Counterexample to Some Conjectures}

\begin{conjecture} 
{\rm (Brady, McCammond, M\"uhlherr, Neumann \cite{B-M-M-N}, Conjecture 8.1)}
Let $W$ be a Coxeter group with Coxeter generating sets $S$ and $S'$. 
If $S^W=(S')^W$, then $\Gamma(W,S)$ is twist equivalent to $\Gamma(W,S')$.
\end{conjecture}

\begin{conjecture} 
{\rm (M\"uhlherr \cite{Muhlherr}, Conjecture 1)} 
Let $W$ be a Coxeter group with Coxeter generating sets $S$ and $S'$ 
such that $S^W=(S')^W$. 
Then there exists a reflection preserving automorphism $\alpha$ of $W$ 
such that $\alpha(S)$ is sharp-angled with respect to $S'$. 
\end{conjecture}

Consider the Coxeter system $(W,S)$ defined by the diagram 
in the left-hand side of Figure 3. 
If we declare the 5-edge $\{b,c\}$ to be bad, 
we have that $\{a,b\}$ is a bad separator of $S$. 
Let $A=\{a,b,c\}$ and let $c' = cabcbac$. 
Then $\beta:\langle A\rangle \to \langle A\rangle$,  
defined by $\beta(a) = b$, $\beta(b) = a$, and $\beta(c) = c'$,  
is an automorphism that is not inner. 
Hence the diagram on the right-hand side of Figure 3 is 
the P-diagram of the system $(W,S')$ obtained from $(W,S)$ 
by the 5-edge angle deformation.  Moreover $S^W = (S')^W$. 
The diagram $\Gamma(W,S)$ supports no elementary twists, 
and so $\Gamma(W,S)$ is not twist equivalent to $\Gamma(W,S')$. 
Thus $(W,S,S')$ is a counterexample to Conjecture 11.1.

There is no reflection preserving automorphism $\alpha$ of $W$ such that $\alpha(S)$ 
is sharp-angled with respect to $S'$, since 
otherwise $\Gamma(W,S)$ would be twist-equivalent to $\Gamma(W,S')$ 
by the Twist Equivalence Theorem. 
Thus $(W,S,S')$ is a counterexample to Conjecture 11.2.
 
\medskip

$$\mbox{
\setlength{\unitlength}{.8cm}
\begin{picture}(16,5)(0,0)
\thicklines
\put(4,1){\circle*{.15}}
\put(4,1){\line(0,1){3}}
\put(4,4){\circle*{.15}}
\put(2,2.5){\line(4,3){2}}
\put(2,2.5){\line(4,-3){2}}
\put(2,2.5){\circle*{.15}}
\put(6,2.5){\line(-4,3){2}}
\put(6,2.5){\line(-4,-3){2}}
\put(6,2.5){\circle*{.15}}
\put(1.4,2.4){$c$}
\put(2.7,1.15){5}
\put(2.7,3.5){3}
\put(4.25,2.3){2}
\put(3.9,.3){$b$}
\put(3.9,4.4){$a$}
\put(5.2,1.15){2}
\put(5.2,3.5){3}
\put(6.4,2.4){$d$}
\put(12,1){\circle*{.15}}
\put(12,1){\line(0,1){3}}
\put(12,4){\circle*{.15}}
\put(10,2.5){\line(4,3){2}}
\put(10,2.5){\line(4,-3){2}}
\put(10,2.5){\circle*{.15}}
\put(14,2.5){\line(-4,3){2}}
\put(14,2.5){\line(-4,-3){2}}
\put(14,2.5){\circle*{.15}}
\put(9.4,2.4){$c'$}
\put(10.7,1.15){5}
\put(10.7,3.5){3}
\put(12.25,2.3){2}
\put(11.9,.3){$a$}
\put(11.9,4.4){$b$}
\put(13.2,1.15){3}
\put(13.2,3.5){2}
\put(14.4,2.4){$d$}
\end{picture}}$$

\medskip

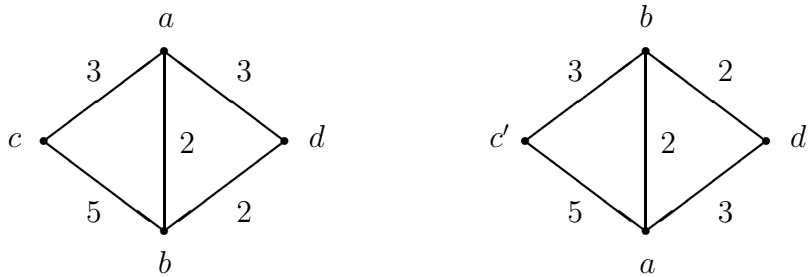
\centerline{\bf Figure 3}

\end{document}